\documentclass[preprint,12pt]{elsarticle}
\usepackage{subfigure}
\usepackage{graphicx,psfrag,epsf}
\usepackage{enumerate}
\usepackage{natbib}
\usepackage{url} 
\usepackage{rotating}
\usepackage[USenglish]{babel} 

\usepackage{amsthm,amsmath,amsfonts,amsthm,bm,booktabs,color}
\usepackage{rotating}
\usepackage{lscape}
\usepackage{graphicx,psfrag,epsf}
\usepackage{float}
\usepackage{natbib}
\usepackage{enumerate}
\usepackage{textcomp}
\usepackage{amsmath}
\usepackage{amsfonts}
\usepackage{amssymb}
\usepackage{latexsym}
\usepackage{verbatim}
\usepackage{times}
\usepackage{color}
\usepackage{bm}
\usepackage{epstopdf}
\usepackage[T1]{fontenc}
\usepackage[utf8]{inputenc}
\usepackage{natbib}

\newtheorem{remark}{Remark}

\graphicspath{{./figuras/}}

\newcommand{\indep}{\;\, \rule[0em]{.03em}{.67em} \hspace{-.25em}
\rule[0em]{.65em}{.03em} \hspace{-.25em}
\rule[0em]{.03em}{.67em}\;\,}

\newcommand{\real}{{\mathbb R}}

\newcommand{\spn}{\mathrm{span}}

\newcommand{\E}{\mathrm{E}}

\newcommand{\tr}{\mathrm{\,tr}}

\newcommand{\I}{\mathbf I}
\newcommand{\Hbf}{\mathbf H}

\DeclareMathOperator*{\argmin}{arg\,min}

\newcommand{\Shat}{\widehat{S}}
\newcommand{\X}{{\mathbf X}}
\newcommand{\x}{{\mathbf x}}

\newcommand{\Abf}{\mathbf A}

\newcommand{\Nbf}{\mathbf N}

\newcommand{\f}{{\mathbf f}}
\newcommand{\fbar}{\bar{\f}}

\newcommand{\z}{{\mathbf z}}
\newcommand{\Z}{{\mathbf Z}}

\newcommand{\F}{{\mathbf F}}
\newcommand{\M}{{\mathbf M}}

\newcommand{\Sbf}{{\mathbf S}}

\newcommand{\A}{{\mathbf A}}

\newcommand{\B}{{\mathbf B}}

\newcommand{\greekbold}[1]{\mbox{\boldmath $#1$}}
\newcommand{\alphabf}{\greekbold{\alpha}}

\newcommand{\alphabfhat}{\widehat{\alphabf}}

\newcommand{\alphabfs}{\greekbold{\scriptstyle \alpha}}

\newcommand{\Thetabfs}{\greekbold{\scriptstyle \Theta}}

\newcommand{\epsilonbf}{\greekbold{\epsilon}}

\newcommand{\Gammabf}{\greekbold{\Gamma}}

\newcommand{\Psibf}{\greekbold{\Psi}}

\newcommand{\Deltabf}{\greekbold{\Delta}}
\newcommand{\Deltabfi}{\greekbold{\Delta}^{-1}}

\newcommand{\Thetabf}{\greekbold{\Theta}}
\newcommand{\Thetabfhat}{\widehat{\greekbold{\Theta}}}
\newcommand{\Deltabfhat}{\widehat{\greekbold{\Delta}}}
\newcommand{\Deltabfs}{{\greekbold{\scriptstyle \Delta}}}

\newcommand{\Omegabf}{\greekbold{\Omega}}

\newcommand{\Sigmabf}{\greekbold{\Sigma}}

\newcommand{\mubf}{\greekbold{\mu}}

\newcommand{\xibf}{\greekbold{\xi}}

\newcommand{\xibfs}{{\greekbold{\scriptstyle \xi}}}
\newcommand{\xibfhat}{{\hat{\xibf}}}

\newcommand{\zetabf}{\greekbold{\zeta}}

\newcommand{\zetabfhat}{{\hat{\zetabf}}}

\newcommand{\res}{\mathrm{res}}
\newcommand{\fit}{\mathrm{fit}}

\newcommand{\Rbb}{{\mathbb R}}

\newcommand{\Psikfy}{\Psibf^{\hspace{-0.05cm}(k-1)} \hspace{-0.05cm} \fbar_{y_i}}
\newcommand{\Omegak}{\Omegabf^{(k-1)}}
\newcommand{\Omegakk}{\Omegabf^{(k)}}
\newcommand{\Deltak}{\Deltabf^{\hspace{-0.1cm}(k-1)}}

\newcommand{\param}{\Thetabf, \Deltabf, \alphabf,  \xibf}
\newcommand{\paramhat}{\Thetabfhat, \Deltabfhat, \alphabfhat, \xibfhat}
\newcommand{\paramehat}{\Deltabfhat, \alphabfhat, \xibfhat}
\newcommand{\parame}{\Deltabf, \alphabf,  \xibf}

\begin{document}

\begin{frontmatter}

\title{Supervised dimension reduction for ordinal predictors}


\author{Liliana Forzani, Pamela Llop, Diego Tomassi \\
Facultad de Ingenier\'ia Qu\'imica, Universidad Nacional del Litoral\\ Researchers of CONICET\\
and \\
Rodrigo Garc\'ia Arancibia \\
Instituto de Econom\'ia Aplicada Litoral (FCE-UNL)}


\bigskip
\begin{abstract}
	\textcolor{black}{
In applications involving ordinal predictors, common approaches to reduce dimensionality are either extensions of unsupervised techniques such as principal component analysis, or variable selection procedures that rely on modeling the regression function.
In this paper, a supervised dimension reduction method tailored to ordered categorical predictors is introduced.
It uses a model-based dimension reduction approach, inspired by extending sufficient dimension reductions to the context of latent Gaussian variables. 
The reduction is chosen without modeling the response as a function of the predictors and does not impose any distributional assumption on the response or on the response given the predictors.
A likelihood-based estimator of the reduction is derived and an iterative expectation-maximization type algorithm is proposed to alleviate the computational load and thus make the method more practical. 
A regularized estimator, which simultaneously achieves variable selection and dimension reduction, is also presented.
Performance of the proposed method is evaluated through simulations and a real data example for socioeconomic index construction, comparing favorably to widespread use techniques.}
\end{abstract}


\end{frontmatter}

\section{Introduction.}
Regression models with ordinal predictors are common in many applications. 
For instance, in economics and the social sciences, ordinal variables are used to predict phenomena like income distribution, poverty, consumption patterns, nutrition, fertility, healthcare decisions, and subjective well-being, among others \cite[e.g.][]{Bollen01,Roy09, Murasko09, Kamakura13, Mazzonna14, Feeny14}. 
In marketing research, customer preferences are used to create automatic recommendation systems, as in the case of Netflix \cite[e.g.][]{bobadilla,roberts}, where the ratings for unseen movies can be predicted using the user's previous ratings and information about the consumer preferences for the whole database. 

In this context, when the number of predictors is large, it is of interest to reduce the dimensionality of the space by combining them into a few variables in order to get efficiency in the estimation as well as an understanding of the model. 
The commonly used dimension reduction techniques for ordinal variables are adaptations of standard principal component analysis (PCA) \cite[][]{Linting09, Kolenikov09}. 
For example, in the case of the {\em \'Indice de Focalizaci\'on de Pobreza} (a socio-economic index commonly used in Latin America), the first normalized principal component is used to predict poverty status, even if this outcome variable was never used to estimate the scaling. 
It is clear, however, that ignoring the response when building such an index can lead to a loss of predictive power compared to the full set of predictors.

A different approach to dimensionality reduction is to perform variable selection on the original set of predictors. 
A method adapted to ordinal predictors is proposed in \cite{Gertheiss10}. 
Despite the fact that this method uses information from the response to achieve variable selection, it performs simultaneously regression modeling by assuming a parametric model for the response as a function of predictors.

For regression and classification tasks it is widely accepted that supervised dimension reduction is a better alternative than PCA-like approaches.
Sufficient dimension reduction (SDR), in particular, has gained interest in recent years as a principled methodology to achieve dimension reduction on the predictors $\X \in \real^p $ without losing information about the response $Y$.
%
%
Formally, for the regression of $Y|\X$, SDR amounts to finding a transformation $R(\X) \in \real^d$, with $d \leq p$, such that the conditional distribution of $Y|\X$ is identical to that of $Y|R(\X)$.
Nevertheless, there is no need to assume a distribution for $Y$ or for $Y|\X$. Thus, the obtained reductions can subsequently be used with any prediction rule. Moreover, when the reduced space has low dimension, it is feasible to plot the response versus the reduced variables. This can play an important role in facilitating model building and understanding \citep{cook96,cook98}.

Most of the methodology in SDR is based on the inverse regression of $\X$ on $Y$, which translates a $p$-dimensional problem of regressing $Y|\X$ into $p$ (easier to model) one-dimensional regressions corresponding to $\X|Y$. 
Estimation in SDR was developed originally for continuous predictors and was based on the moments of the conditional distribution of $\X|Y$ (SIR, \cite{li91}; SAVE, \cite{cookweisberg}; pHd, \cite{li92};  PIR \cite{bura01}; MAVE, \cite{Xia02,Li05,CookNi,Zhu06,CookLi}; DR, \cite{LiWang}, see also \cite{CookLee,CookYin,chiaro,CookNi,yin-li-cook} and \cite{cook98}, where much of its terminology was introduced). 
Later, \cite{cook_2007} introduced the so called  \textit{model-based} inverse regression of $\X|Y$ (see also \cite{cook_forzani_2008,cook_forzani_2009}).
The main advantage of this approach is that provides an estimator of the sufficient reduction that contains all the information in $\X$ that is relevant to $Y$, allowing maximum likelihood estimators which are optimal in terms of efficiency and $\sqrt{n}$-consistent under mild conditions when the model holds. 

Along the lines of the model-based SDR approach, the up to date methodology is for predictors belonging to a general exponential family of distributions (See \cite{BDF}).
Then, when attempting to apply SDR to ordinal predictors, a first approach could be to treat them as polythomic variables, ignoring their natural order. Then a multinomial distribution can be postulated over them, which can be treated as member of the exponential family.
However, ordered variables usually do not follow a multinomial distribution and the order information is lost when treating them as multinomial. 
There have been attempts to use dummy variables to deal with ordinal data, but this procedure has been shown to introduce spurious correlations \citep{Kolenikov09}.
Another approach is to treat the ordered predictors as a discretization of some underlying continuous random variable. 
This technique is commonly used in the social sciences and is known as the \textit{latent variable model}. 
In this context, the latent variables are usually modeled as normally or logistically distributed, obtaining the so-called ordered probit and logit models, respectively \citep{greene_2010, long_1997}. 
While for each scientific phenomenon the latent variable can take a particular meaning (e.g., utility in economic choice problems, liability to diseases in genetics, or tolerance to a drug in toxicology), a general interpretation of a latent variable  may be the {\it propensity} to observe a certain value $j$ of an ordered categorical variable \citep{Skrondal04}. 
Regardless of their philosophical meaning and the criticisms about their real existence, latent variables are very useful for generating distributions for modeling, hence their widespread use.

In this paper, we develop a supervised dimension reduction method for ordinal predictors, based on the SDR for the regression of the response given the underlying normal latent variables. 
Under this context, we present a maximum likelihood estimator of the reduction and we propose an approximate expectation maximization (EM) algorithm for its practical computation, which is 
close to recent developments in graphical models for ordinal data \citep{levina_2014} and allows for computationally efficient estimation without losing accuracy.

The rest of this paper is organized as follows. 
In Section~\ref{model} we describe the inverse regression model for ordinal data and its dimensionality reduction. 
In Section~\ref{estimation} we derive the Maximum likelihood estimates of the reduction and we also present a variable selection method. Section~\ref{dimension} is dedicated to developing a permutation test for choosing the dimension for the reduction.
Simulation results are presented in Section~\ref{simulation}.
Section~\ref{realdata} contains a socio-economic application using the methodology developed in this paper to create a socio-economic status (SES) index from ordinal predictors.
Finally, a concluding discussion is given in  Section~\ref{Conclusion}.
All proofs and other supporting material are given in the appendices.
Matlab codes for the algorithm and simulations are available at {\color{black}{\small{\texttt{http://www.fiq.unl.edu.ar/pages/investigacion/\ investigacion-reproducible/grupo-de-estadistica-statistics-group.php}}}}.

\section{Model} \label{model}
Let us consider the regression of a response $Y \in \Rbb$ on a predictor {$\X= (X_1,X_2,\dots,X_p)^T$}, where each $X_j$, $j =1,\ldots,p$ is an ordered categorical variable, i.e., $X_j \in \{1,\dots,{G_j}\}$, $j =1,\ldots,p$. To state a dimension reduction of $\X$ inspired by the model-based SDR approach (see \cite{cook_2007}), we should model the inverse regression of $\X$ on $Y$. However, as we stated in the introduction, the model-based SDR techniques deal with continuous predictors. Therefore, in order to frame our problem in that context, we will assume the existence of a $p$-dimensional vector of unobserved underlying continuous latent variables $\Z = (Z_1,Z_2,\dots,Z_p)^T$, with $E(\Z) =0$, such that each {\color{black}{observed $X_j$ is 
a discretizing of $Z_j$}} as follows.
{\color{black}{There exists 
 a set of thresholds $\Thetabf^{(j)} = \{\theta^{(j)}_0, \theta^{(j)}_1,\dots,\theta^{(j)}_{G_j}\}$, that split the real line in disjoints  intervals $-\infty=\theta^{(j)}_0 < \theta^{(j)}_1 < \dots < \theta^{(j)}_{G_{j-1}} < \theta^{(j)}_{G_j}=+\infty$ and }} 
\begin{equation}\label{eq:XfromZ}
X_j = \sum_{g=1}^{G_j} g \mathbb{I}(\theta^{(j)}_{g-1}\le Z_j < \theta^{(j)}_g),
\end{equation}
where $\mathbb I(A)$ is the indicator function of the set $A$. Therefore, $X_j = g \Leftrightarrow Z_j \in [\theta^{(j)}_{g-1},\theta^{(j)}_{g})$ and $P(X_j = g) = P(\theta^{(j)}_{g-1}\le Z_j < \theta^{(j)}_g).$

In the framework of model-based inverse regression, we adopt, following \cite{cook_forzani_2008}
{\color{black}{that the variable $\Z$ given $Y$ is normal with mean depending on $Y$ and constant variance, i.e.}}
\begin{equation}\label{model1}
\Z|Y = \mubf_Y+ \epsilonbf, 
\end{equation}
where $\mubf_Y = E(\Z|Y)$ and the error $\epsilonbf$ is independent of $Y$, normally distributed with mean $\bm 0$ and covariance (positive definite) matrix $\Deltabf$. 
As usual in latent variable models for ordinal data (see \cite{jackman_2009}), additionally to $E(\Z) =0$, we set the diagonal $[\Deltabf]_{jj} \doteq \delta_j = 1$ in order to allow for model identification.

{\color{black}{
Since $E(\Z|Y)$ depends on $Y$ we could model that dependence as a function of  $\f_Y \in \mathbb R^r$  a  vector of $r$ known functions with $E((\f_Y-E(\f_Y))( \f_Y- E(\f_Y))^T)$.
Under this model, each coordinate of $\Z|Y$ follows a linear model with predictor vector $\f_Y$ and therefore, when $Y$ is quantitative, we can use inverse plots to get information about the choice of $\f_y$, which is not possible in the regression of $Y$ on $\X$. When $Y$ is continuous, $\f_y$ usually will be a flexible set of basis functions, like polynomial terms in
$Y$, which may also be used when it is impractical to apply graphical methods to
all of the predictors. 
When $Y$ is categorical and takes values $\{C_1,\ldots,C_h\}$, we can set $r=h-1$ and specify the $j$th element of $\f_y$ to be $\mathbb I(y \in C_j)$, $j=1,\ldots,h$. When $Y$ is continuous, we can also {\it slice} its values into $h$ categories $\{C_1,\ldots,C_h\}$ and then specify the $j$th coordinate of $\f_y$ as for the case of a categorical $Y$. 
For more details see \cite{Adragni2009}. As a consequence, model (\ref{model1}) can be expressed as 
\begin{eqnarray}\label{model22}
\Z|Y  =   \Gammabf  \{\f_Y-E(\f_Y)\} + \epsilonbf,
\end{eqnarray}
where $\epsilonbf$ is independent of $Y$, normally distributed with mean $\bm 0$ and covariance (positive definite) matrix $\Deltabf$. }}

{\color{black} For the regression of $Y$ on the continuous latent variable $\Z$, under model (\ref{model22}) the minimal SDR is $R(\Z)= \alphabf^T \Z$, with $\alphabf $ a basis for $ \Deltabf^{-1} \spn (\Gammabf)$ by \cite{cook_forzani_2008}, Theorem 2.1. Note that if $R(\Z) = \alphabf^T  \Z$ is a sufficient reduction, then $R(\Z) =\A \alphabf^T \Z$ is a sufficient reduction for any invertible $\A \in {\mathbb R}^{d\times d}$ \citep{cook98}. 
Therefore what is identifiable is the span of $  \alphabf$, not $\alphabf$ itself. In the SDR literature, the identifiable parameter $\hbox{span}(\alphabf)$ is called a {\it sufficient reduction subspace}. If dim$(\spn (\Gammabf))=d\le \min \{r, p\}$,
 (\ref{model22}) can be re-written as
\begin{eqnarray}\label{model2}
\Z|Y  =  \Deltabf \alphabf \xibf\{\f_Y-E(\f_Y)\} + \epsilonbf,
\end{eqnarray}
where  $\alphabf \in \Rbb^{p\times d}$ with $d \leq p$ is a semi-orthogonal matrix whose columns form a basis for the 
$d$-dimensional subspace $\Deltabf^{-1} \spn (\Gammabf)$,
$\xibf \in\Rbb^{d\times r}$ is a full rank $d$ matrix with $r\geq d$ (see \cite{cook_forzani_2008}, \cite{Adragni2009}).

Coming back to our problem of interest, in order to propose a supervised dimension reduction for the regression of $Y|\X$ let us observe that, since $\X$ is a function of $\Z$, $\alphabf^T  \Z$ will be also the sufficient dimension reduction for $Y|\X$, i.e. $Y \indep \X | \alphabf^T  \Z$ (see Proposition 4.5 in \cite{cook98}). {\textcolor{black}{However, since $\Z$ is unobservable, {\color{black}{and the only information available is $\X$, we take the conditional expectation of $\alphabf^T \Z$ given $\X$, instead of $\alphabf^T \Z$
for the reduction of $Y|\X$ since it is the best predictor of $\alphabf^T\Z$ in terms of minimum Mean Square Error. Therefore, for the regression of $Y$ on $X$, the proposed supervised dimension reduction will be}}
\begin{equation}\label{nuevared}
R(\X)=  E(\alphabf^T\Z|\X).
\end{equation}}

\textcolor{black}{
\begin{remark}
Observe that in this case, regardless of the encoding of $\X$, the reduction $R(\X)$ is completely identified since for each $j =1,\ldots,p$, $E(Z_j)=0$ and $X_j = g \Leftrightarrow Z_j \in [\theta^{(j)}_{g-1},\theta^{(j)}_{g})$ and therefore, whatever the coding of $X_j$ is, the underlying (and as a consequence the thresholds) does not change.
\end{remark}}

\begin{figure}
	\centering
	\includegraphics[scale=.8]{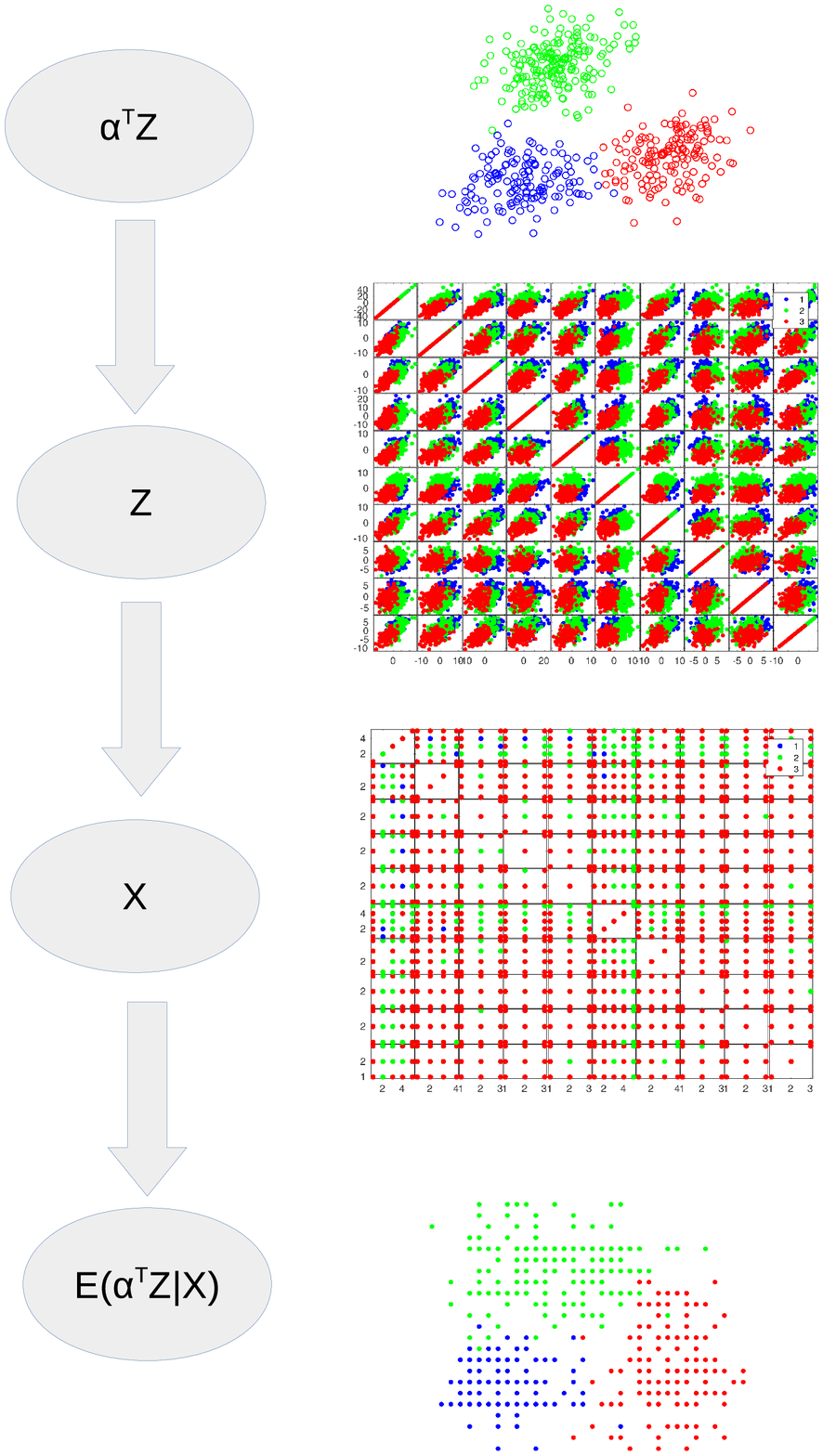}
	\caption{Data generation and processing chain according to the assumed model and proposed dimension reduction scheme.}
	\label{fig:ilu}
\end{figure}

Figure \ref{fig:ilu} helps to understand how the proposed method works, showing an example corresponding to a categorical outcome with three nominal values $Y\in\{1,2,3\}$. Suppose we are able to observe the underlying continuous variables $\Z$ and that their distribution follows model (\ref{model2}), with $\alphabf \in \real^{p\times 2}$. It means that the characteristic information needed to discriminate between the three groups lies actually in the two-dimensional subspace spanned by the columns of $\alphabf$. Top panel of Figure~\ref{fig:ilu} shows such information, plotting the coordinates of $\alphabf^T\Z$. For each group indexed by $Y$, the data fall in a cluster well separated from the others.
If we are not allowed to observe this reduced subspace but the complete underlying predictors $\Z$, we would have a situation as described by the scatter plots between pairs of predictors depicted in the second panel of Figure \ref{fig:ilu}. Despite we can still see some separation between clusters, it is not as clear as in the sufficient low-dimensional subspace.
In real scenarios with ordinal data, according to the assumed model we do not have access to observe $\Z$ either, but a discretized version $\X$ which is a function of the underlying $\Z$ throught the set of fixed thresholds $\Thetabf$. This situation is illustrated in the third panel. It is clear from the figure that the continuous values of $\Z$ collapse into a discrete set of values in $\X$ and it is now much harder to discriminate between the groups indexed by $Y$.
Nevertheless, the dimension reduction approach proposed in this paper projects the data again onto a 2-dimensional subspace, as shown in panel at the bottom of Figure~\ref{fig:ilu}. Note that clear separation between clusters is recovered; indeed, the information available in this subspace closely resembles that in the characteristic subspace spanned by $\alphabf$ \textcolor{black}{(compare the first and last panels)}.

\begin{figure}
	\centering
	\includegraphics[scale=.5]{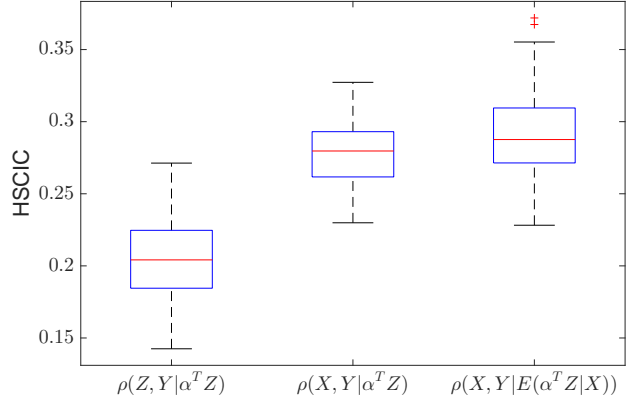}
	\caption{Residual dependence between the response and the predictors, given the reduction. Hilbert-Schmidt condicional independence criterion (HSCIC) is used as measure of conditional dependence.  }
	\label{fig:simu4dcov}
\end{figure}

It is interesting to see also how well the proposed reduction captures the information about $Y$ that is available in $\X$. We know that for the underlying variables, $\Z \indep Y \mid \alphabf^T\Z$. If we can measure the residual dependence between $\Z$ and $Y$ given $\alphabf^T\Z$ using a \textcolor{black}{suitable measure of statistical independence} $\rho(\Z,Y\mid\alphabf^T\Z)$ , we will have $\rho(\Z,Y\mid\alphabf^T\Z)=0$ in the population. 
In practice, for a finite random sample \textcolor{black}{ of $(Y,\Z,\alphabf^T\Z)$}, this value will be greater than cero. This is shown in the left-most boxplot in Figure \ref{fig:simu4dcov}. The boxplot was computed using 100 simulated data sets and choosing for the measure $\rho$ the Hilbert-Schmidt condicional independence criterion introduced in \textcolor{black}{\cite{hscic}}. Isotropic Gaussian kernels are used to embed the observations into a RKHS. Kernel bandwiths are set to the median of the pairwise distances between samples.
Since $\alphabf^T\Z$ is still a sufficient reduction for the regression $Y|\X$, in the population we would also have $\rho(\X,Y\mid \alphabf^T\Z)=0$. A sample estimate of this quantity, using the same data and dependence measure as before, is shown in the second boxplot of the figure.
Finally, \textcolor{black}{if we compute the same empirical measure, with the same random sample but for the practical reduction $R(\X)=E(\alphabf^T\Z\mid\X)$ instead of $\alphabf^T\Z$, we obtain the boxplot shown on the right of the figure. This figure shows that on average the empricial version of $\rho(\X,Y\mid E(\alphabf^T\Z|\X))$ is a little bigger than the empirical version of $\rho(\X,Y\mid \alphabf^T\Z)$ for the same random sample of $(Y,\Z,\alphabf^T\Z)$, albeit they are very close.
This suggests that, even if we cannot claim sufficiency of the proposed reduction, it is really close to the ideal unattainable reduction $\alpha^T\Z$.}

\section{Estimation} \label{estimation}
\textcolor{black}{For the supervised dimension reduction given in (\ref{nuevared}), we need to estimate the semiorthogonal basis matrix $\alphabf$.} If $\Z$ were observed, the maximum likelihood estimator of $\alphabf$
would be the one derived in \cite{cook_forzani_2008}. 
That is, $\alphabfhat =  \widetilde{\Sigmabf}^{-1/2} {\mathbf v}$, where ${\mathbf v}$ are the first $d$ eigenvectors of the symmetric matrix $\widetilde{\Sigmabf}^{-1/2}\widetilde{\Sigmabf}_{\fit} \widetilde{\Sigmabf}^{-1/2}$, $\widetilde{\Sigmabf} $ is the sample marginal covariance of the predictors and
$\widetilde{\Sigmabf}_{\fit} $ is the sample covariance of the fitted values of the regression $\Z|\f_y$. This estimation procedure is called Principal Fitted Components (PFC). However, $\Z$ is not observed, and therefore
the sample covariance matrices (marginal and fitted) cannot be estimated directly. In view of the robustness proven in \cite{cook_forzani_2008}, we could consider applying the methodology directly to $\X$ in a naive way and it still will obtain a $\sqrt{n}$ consistent estimator. This approach will be the initial value of our algorithm to obtain the maximum likelihood estimate under the true model. 
%

For the estimation, let as assume we have a random sample of $n$ points $(y_i,\x_i)$ drawn from the joint distibution of $(Y,\X)$ following model  (\ref{model2})  
and that the dimension $d$ of the reduction is known (later in Section \ref{dimension} we will consider how to infer it).
\textcolor{black}{In what follows, we will call $\Thetabf \doteq \{\Thetabf^{(1)},\ldots,\Thetabf^{(p)} \} =\{\theta^{(1)}_0,\dots,\theta^{(1)}_{G_1},\dots \ldots, \theta^{(p)}_0,\dots,\theta^{(p)}_{G_p}\}$ and $C(\X,\Thetabf) = [\theta^{(1)}_{X_1-1},\theta^{(1)}_{X_1})\times \dots \times [\theta^{(p)}_{X_p-1},\theta^{(p)}_{X_p})$.} In order to obtain the estimator we need 
to maximize the log-likelihood function of the observed data
\begin{equation}\label{max}
\sum_{i=1}^n \log f_{\X} (\x_i|y_i; \param).
\end{equation}
%
Since $\Z|Y$ is normally distributed, using (\ref{eq:XfromZ}) we can compute, for 
each $i$, the truncated unnormalized density $f_{\X,\Z}(\x_i, \z_i|y_i;\param)$ as
\[ 
f_{\X,\Z}(\x_i,\z_i|y_i;\param) = (2\pi)^{-p/2} |\Deltabf|^{-1/2} e^{-\frac{1}{2}  \tr(\Deltabfs^{-1} (\z_i-\Deltabfs\alphabfs \xibfs \fbar_{y_i})(\z_i-\Deltabfs \alphabfs \xibfs \fbar_{y_i})^T)}I_{\{\z_i \in C(\x_i,\Thetabfs)\}},
\]
where $\fbar_{y_i} \doteq \f_{y_i} - n^{-1}\sum_{i=1}^n \f_{y_i}$. Therefore, for each $i$, the unnormalized marginal density $\X|Y$ will be 
\begin{equation*}\label{integral}
f_{\X}(\x_i|y_i;\param) =  \int \, f_{\X,\Z}(\x_i,\z_i|y_i;\param) \, d\z_i.
\end{equation*}
As an exact computation of the likelihood is difficult due to the multiple integrals involved, maximum likelihood estimates are often obtained using an iterative expectation-maximization (EM) algorithm.
This is a common choice for models with latent variables, since it exploits the reduced complexity of computing the joint likelihood of the complete data $(\X,\Z)$. We will follow this approach in the present paper. The corresponding algorithm is described below.

\subsection{Algorithm} \label{algorithm}
In this section we present the EM algorithm, closely related to the one given in \cite{levina_2014}, to estimate the parameters in model (\ref{model2}). 
Throughout this section, we will use superscripts $A^{(k)}$ to indicate the value of quantity $A$ at the $k$th iteration of the algorithm.
In addition, to make notation easier, let us collect the parameters \textcolor{black}{$\Deltabf,\alphabf,\xibf$} into a single parameter vector \textcolor{black}{$\Omegabf\doteq\{\Deltabf,\alphabf,\xibf\}$}.
The procedure starts with \underline{Step 0}, where we initialize $\Omegabf^{(0)}$ using the estimators obtained from \texttt{PFC} applied to $\X$. 
Then, the algorithm iterates between the following two steps until convergence is reached: \underline{Step 1} is devoted to estimating $\Thetabf^{(k)}$ given 
$\Omegabf^{(k-1)}$ and
\underline{Step 2} to getting $\Omegabf^{(k)}$ by maximizing the
conditional expectation (given $\Omegabf^{(k-1)}$ and $\Thetabf^{(k)}$) of the joint log-likelihood  (\ref{max}). 
This step is properly the EM step.

\textbf{\underline{Step 1:} Estimation of $\Thetabf$:} Given $\Omegabf^{(k-1)}$  from \underline{Step 
0} or from a previous iteration, let $\Psibf^{(k-1)} 
\doteq \Deltabf^{(k-1)}\alphabf^{(k-1)}\xibf^{(k-1)}$.
For each $j=1,\dots, p$ and $g=1, \dots, G_j$ define
\[ 
L^{(j)}_g(\theta) \doteq \# \{ i: x_{ij} \le g\} - \sum_{i=1}^n \Phi \left(\frac{\theta- 
\Psibf^{(k-1)}_j \fbar_{y_i}}{\delta^{(k-1)}_j} \right),
\]
where $\Phi$ is the cumulative distribution function of the standard normal, for each $j$, $\delta^{(k-1)}_j=(\Deltabf^{(k-1)})_{jj}$, $\Psibf^{(k-1)}_j$ indicates 
the $j$th row  of $\Psibf^{(k-1)}$,
$x_{ij} $ is the $j$th coordinate of $\x_i$,
and $\# S$ indicates the cardinality of the set $S$.
Then, take $\hat \theta^{(j)}_0 = -\infty=$,
$\hat \theta^{(j)}_{G_j}=+\infty$.
For $g = 1, \ldots, G_j -1$, assign to $\hat{\theta}^{(j)}_g$ the
unique solution of the equation $L^{(j)}_g(\theta) = 0$.
Set $\Thetabf^{(k)} = \{\theta^{(1)}_0,\dots,\theta^{(1)}_{G_1},\dots,
\theta^{(p)}_0,\dots,\theta^{(p)}_{G_p}\}$. 
Here the definition of $L^{(j)}_g$ is based on the normality assumption on the conditional distribution of the underlying continuous variable. 
More precisely, we define $L^{(j)}_g$ as a {\it search} function of thresholds using the underlying (normal) cumulative distribution function.

\textbf{\underline{Step 2:} Estimation of  $\Deltabf^{(k)}, \alphabf^{(k)},  \xibf^{(k)}$ :} Given $\Thetabf^{(k)}$ computed in \underline{Step 1} and $\Omegak$  from \underline{Step 0} or from a previous iteration, we apply the EM algorithm to maximize (\ref{max}). 
The EM algorithm consist in finding  $\Omegakk$ that maximize over $\Omegabf$
\begin{equation}\label{QQQ1}
Q(\Omegabf | \Omegak) = \sum_{i=1}^n E_{\z_i |y_i;\tiny{\Omegak}} \left[ \log f_{\x_i,\z_i}(\x_i,\z_i|y_i;\Omegabf)\big |y_i;\Omegak \right].
\end{equation}
These produce
\begin{eqnarray*}\label{alfahat}
\alphabf^{(k)} &=& \Sbf^{-1/2} \zetabfhat_d \Nbf,\\
(\Deltabf^{-1})^{(k)} &=&\Sbf^{-1} + \alphabf^{(k)} ((\alphabf^{(k)})^{T} \Sbf_{\res} \alphabf^{(k)})^{-1} (\alphabf^{(k)})^{T}  - \alphabf^{(k)} ((\alphabf^{(k)})^{T} \Sbf \alphabf^{(k)})^{-1} (\alphabf^{(k)})^{T}
,\\
\xibf^{(k)}& =& ( (\alphabf^{(k)})^T  \Deltabf^{(k)} \alphabf^{(k)})^{-1}( \alphabf^{(k)})^T \M^T \F (\F^T 
\F)^{-1},
\end{eqnarray*}
where $\Nbf$ is a matrix such that $(\alphabf^{(k)})^T \alphabf^{(k)}= 
\mathbf I_d$ and the $\zetabfhat_d$ are the first $d$ eigenvectors of 
$ \Sbf^{-1/2}\Sbf_{\fit} \Sbf^{-1/2}$, where the matrices $\Sbf \in \Rbb^{p\times 
p}$ and $\Sbf_{\fit} \in \Rbb^{p\times  p}$ are given by 
\[
\Sbf = \frac{1}{n} \sum_{i=1}^n E_{\z_i| \x_i,y_i;\tiny\Omegak} (\z_i \z_i^T|\x_i, y_i;\Omegak) \quad \text{and} \quad \Sbf_{\fit} =  n^{-1}\M^T \F (\F^T \F)^{-1} \F^T \M 
\]
with $\F \in \Rbb^{n\times r}$ and $\M\in \Rbb^{n\times p}$ matrices whose transposes are given by $\F^T = [\fbar_{y_1}, \ldots, \fbar_{y_n}]$ and $\M^T = [E_{\z_1|\x_1,y_1;\tiny\Omegak} (\z_1|\x_1, y_1;\Omegak), \ldots, E_{\z_n|\x_n,y_n;\tiny\Omegak} (\z_n|\x_n, y_n;\Omegak))]$ and the residual matrix $\Sbf_{\res}$ is defined by $\Sbf_{\res}=\Sbf-\Sbf_{\fit}$. Details of the derivation and the EM algorithm are given in \ref{appendixA}.

\textbf{\underline{Step 3:}} Check convergence. 
If it is not reached, go to \textbf{Step 1}. 
We check convergence simply by looking to see whether $Q(\Omegabf^{(k)}|\Omegabf^{(k-1)})$ stops increasing from one iteration to the next. 
Specifically, we check whether $(Q(\Omegabf^{(k)}|\Omegabf^{(k-1)}) - Q(\Omegabf^{(k-1)}|\Omegabf^{(k-2)}))/Q(\Omegabf^{(k-1)}|\Omegabf^{(k-2)}) < \epsilon$, with $\epsilon$ typically set to $10^{-6}$.

\subsection{Estimation with variable selection} \label{sec:lasso}

\textcolor{black}{When we compute the linear 
combinations implied in (\ref{nuevared}), we need to include all the original variables.
This means that even non-relevant or redundant variables are included in the final model, making it harder to interpret.}
To overcome this limitation, we can perform variable 
selection, in this way obtaining linear combinations that include only the active, relevant variables. 
Following \cite{chen_2010}, the maximization of (\ref{QQQ1}) 
is equivalent to finding, in each iteration, 
\begin{equation} \label{alphavs}
\alphabf^{(k)} = \argmin_{\alphabf} \big\{-\tr (\alphabf^T \Sbf_{\fit}\alphabf)\big\}, 
\hspace{1cm} 
\text{subject to } \alphabf^T \Sbf\alphabf = \mathbf I_{d}. 
\end{equation}
To induce variable selection in dimension reduction, we introduce a group-lasso type penalty since, in order not to choose a particular variable $X_j$, we need to make the whole $j$th row of $\alphabf$, $\alphabf_j$, equal to 0. 
For that, following
\cite{chen_2010}, we use a mixed $\ell_1/\ell_2$ norm, where the inner norm is the $\ell_2$ norm of each row of $\alphabf$. 
Adding the penalty term to (\ref{alphavs}), we get
\[
\alphabf^{(k)} =  \argmin_{\alphabf} \Big\{-\tr (\alphabf^T \Sbf^{-1/2}\Sbf_{\fit} \Sbf^{-1/2} 
\alphabf)+ \lambda \sum_{i=1}^p ||\alphabf_i||_2\Big\}, \hspace{1cm} \text{subject to } \alphabf^T \Sbf\alphabf = \mathbf I_{d}.
\]

The parameter $\lambda$ can be found using an information
criterion, such as Akaike's (AIC) or Bayes' (BIC) criteria. 
Details are given in \cite{chen_2010}. 
Another approach is to find the value $\lambda^{*}$ that minimizes the prediction error via a cross-validation experiment, but this requires adopting a specific prediction rule.

It is interesting to note that this procedure performs at the same time variable selection and dimension reduction without modeling the regression for $Y|E(\Z|\X)$ or $Y|\X$. 
Thus, the obtained reduction can be used later with any prediction rule of choice. 
This is different, for instance, from the approach proposed in \cite{Gertheiss10}, where the variable selection is driven by a particular regression model.

\subsection{Computing the reduction}

In order to compute the reduction (\ref{nuevared}), we need to estimate  $\alphabf$ and $E(\Z|\X)$. In the preceding paragraph we focused on computing an estimate of $\alphabf$. To estimate $E(\Z|\X)$, observe that, when the response $Y$ is discrete, using Bayes' rule we get
\begin{align*}
E(\Z|\X) &= E_{Y|\X}\big(E(\Z|\X,Y=y)\big)\\ &= \sum_{y \in S_Y} P(Y=y|\X)E(\Z|\X,y)
\\ &= \sum_{y \in S_Y} w(y)E(\Z|\X,y),
\end{align*}
where, 
\[
w(y) = \frac{P(\X|Y=y) P(Y=y)}{\sum_{y \in S_Y}P(\X|Y=y) P(Y=y)},
\]
Then, to estimate $E(\Z|\X)$ we take: $\widehat{E}(\Z|\X,y)=\M$ from the Step 2 of the EM algoritm (Section \ref{algorithm}); $\widehat{P}(\X|Y=y) = \int_{C(\X, \widehat{\Thetabf})} f_{\Z|Y}(\z|y;\paramehat) \, d\z$ with $\paramhat$ obtained from Steps 1 and 2 of Section \ref{algorithm} or \ref{sec:lasso}; and $\widehat{P}(Y=y)$ from the sample. When the response is continuous, we can simply slice $Y$ in $h$ bins and use the previous procedure. 
Note also that the sample space of $\X$ is finite and so the sample space of $E(\alphabf^T\Z|\X)$ is finite too. Thus, a priori we can tabulate those values for future use and avoid computations when we need to reduce a new instance of $\X$. Nevertheless, when $p$ is moderately large and there are several ordered categories for each variable $X_j$, the amount of memory needed to store such a look up table can become too large in practice. For instance, if $p=20$ and we have $G_j=3$ for each $X_j$, to store all the values of $E(\alphabf^T\Z|\X)$ in double precision we need around 26 GB of memory.

\section{Choosing the dimension} \label{dimension}
Our developments in Section~\ref{estimation} assumed that the dimension $d \le \min (p,r)$ of the reduction subspace was known. In practical settings, this dimension should be inferred from the data. For model-based SDR, which allows for likelihood computation, likelihood-ratio and information criteria such as AIC or BIC have been proposed to drive the selection of $d$ \citep{cook_forzani_2008}. The accuracy of these methods, however, is not robust to deviations from the assumed model. When the main goal is prediction, a common choice is to assess different values of $d$ according to their performance at predicting out-of-sample cases in a cross validation setting. 
The value of $d$ picked is the one that achieves the minimum prediction error.

Another way to choose the dimension is  via permutation tests, as introduced in \cite{CookYin}.
For that, assume that $\alphabf \in \real^{p\times m}$ and $(\alphabf , \alphabf_0) \in \real^{p \times p}$ is unitary. 
A permutation  test relies on the fact that $\mathcal{S}_{\alphabfs}$ is a sufficient dimension reduction subspace for the regression of $Y$ on $\Z$ whenever $(Y,\alphabf^T\Z) \indep \alphabf_0^T\Z$. Note that this implies $m \geq d$. 
For this test, we consider the statistic $\hat{\Lambda}_m = 2({\mathcal Q}_p(\I_p) - \mathcal{Q}_m({\alphabfhat}))$, where ${\mathcal Q}_{r}$ is given by the $Q$ function in (\ref{QQQ1}), evaluated at the estimator obtained in the EM algorithm given above, for a fixed dimension of $\mathcal{S}_{\alphabfs}$. Set $m=0$. A procedure adapted to ordinal data to infer $d$ via permutation testing involves the following steps:
\begin{enumerate}[(i)]
\item Obtain $\alphabfhat$, the MLE of $\alphabf$  and compute $\hat{\Lambda}_m$. 
Obtain also $\alphabfhat_0$.
\item For the data $(Y_i,E(\alphabfhat^T\Z|\X_i),E(\alphabfhat_0^T\Z|\X_i))$, permute the columns corresponding to $E(\alphabfhat_0^T\Z|\X_i)$ to get a new sample $(Y_i,E(\Z|\X_i)^*)$. 
For the new data, obtain the MLE and compute $\hat{\Lambda}_m^*$.
\item Repeat step 2, $B$ times. 

\item Compute the fraction of the $\hat{\Lambda}_m^*$ that exceed $\hat{\Lambda}_m$. 
If this value is smaller than the chosen significance level and if $m < \min (r,p)$ set $m=m+1$ and go to step 1. 
Otherwise, return $d=m$.
\end{enumerate}

In this way, the inferred $d$ is the smallest $m$ that fails to reject the null hypothesis of independence between $(Y,E(\alphabfhat^T\Z|\X_i))$ and $E(\alphabfhat_0^T\Z|\X_i)$.

\section{Simulations} \label{simulation}

In this section we illustrate the performance of the proposed method using simulated data.
A critical aspect of the implementation is the computation of the E-step of the EM algorithm. Exact computation of the truncated moments involved in this step would make the proposed method infeasible in practice even for a dimensionality of the predictors of order $5 < p < 10$, depending on the number of ordered categories.
To address this, we implement an approximate estimation method adapted from \cite{levina_2014}. 
The main idea is to use a recursion to iteratively compute the truncated moments of a 
multivariate normal distribution. The derivation and further details are given in  
\ref{appendiceC}. The first step is to validate the proposed approximate method on the E-step by comparing its 
performance with the exact computation of the truncated moments.
Then, we compare the performance of the proposed SDR\textcolor{black}{-based} method against  standard methodology 
developed for continuous data. Taking advantage of the computational savings obtained with the approximate E-step computation, we then illustrate the performance of the proposed stra\-tegy to infer  the dimension of the dimension reduction subspace using permutation testing and cross validation.
Finally, we illustrate the performance of the regularized estimator proposed in 
Section~\ref{sec:lasso} in a prediction task.

\subsection{Validation of the proposed algorithm for the E-step of the algorithm}\label{sec:validation}
The most demanding part of the proposed method is the computation of the truncated moments of the multivariate normal distributions in the E-step. 
The approximate iterative method proposed for its computation is a main ingredient to allow the application of the methodology in practical settings. 
In this section we validate this strategy by comparing the approximate computation against the exact computation using the algorithm proposed in \cite{LeeScott}.  Since the exact computation involves a high computational cost even for moderate dimensions of the predictor vectors, we set $p=5$ and $n=100$. In addition, we set $G_j = 4$ for $j=1,2,\dots,5$. The data was generated according to \eqref{model2}, with $Y\sim N(0,1)$. 
For the basis matrix $\alphabf$, we set
$\sqrt{p}\;\alphabf = (\mathbf{1}_p \;\; \hbox{sign}(\mathbf{e}))$,
with $\mathbf{1}_p$ a column vector of ones of size $p=5$ and
$\mathbf{e}\sim N_p(\mathbf{0},\I)$.
For the covariance matrix $\Deltabf$, we set
$\Deltabf = \I + \alphabf\B\alphabf^T$, with
$\B$ a $2\times 2$ symmetric random matrix fixed at the outset. We also  chose a polynomial basis for $\f_Y$, with $r=2$. 
The same choice of $\f_Y$ was used for the estimation. The experiment was replicated 100 times. 
In each run, the same training sample was used with both methods.

To assess the accuracy of the estimation,  we measured the angle between the subspace spanned by the true $\alphabf$ and the one spanned by the estimate $\alphabfhat$. This quantity ranges from 0 degrees if the two subspaces are identical to 90 degrees if they do not share any information. The average angle obtained with the exact method was $10.21$ degrees with a standard deviation of $6.37$ degrees, whereas for the estimate obtained with the approximate method the average angle was $12.65$ degrees with a standard deviation of $5.70$ degrees. 
The 95\% confidence interval for the average difference between the angles obtained with both methods is $(2.07,2.81)$ degrees. These values suggest that the price to pay for the introduction of the approximate computation is very small. 

It is also illustrative to see the impact of the approximation on prediction. 
Using plain linear regression for $Y|E(\alphabfhat^T\Z|\X)$, the MSE of the residuals averaged over the 100 runs is $0.863$ when estimating the reduction using the approximate E-step, whereas it is $0.811$ when using the exact method. 
In both cases the standard deviation of the averaged MSE is $0.04$. 
Thus, the differense in the MSE obtained with the approximate method represents less than 1.5\% of the average MSE obtained with the exact method. 
The importance of the approximate method for the E-step is better understood when noting the big difference in computing time. 
Using plain MATLAB implementations for both methods, computation with the exact method takes on average $1.13\times 10^2$ seconds for each run, while with the approximate method this time was reduced to $0.26$ seconds, a difference of three orders of magnitude. 
Overall, these results show that the approximate method to compute the truncated moments is a viable alternative: it reduces the computing time in practical applications without a significant loss in accuracy.

\subsection{Performance of the proposed method} \label{sec:simu1}
\textit{}
In this section we assess the performance of the proposed method on simulated data.
For this example, we set $p=20$, $d=2$ and a polynomial basis with $r=2$ for $\f_Y$. As in Section~\ref{sec:validation}, we generate the data according to \eqref{model2}, with $Y\sim N(0,1)$, $\sqrt{p}\,\alphabf = (\mathbf{1}_p \;\; \hbox{sign}(\mathbf{e}))$, with $\mathbf{e}\sim N_p(\mathbf{0},\I)$, and covariance matrix $\Deltabf = \I + \alphabf\B\alphabf^T$, with $\B$ a $d\times d$ symmetric random matrix fixed at the outset. The values of $G_j$ in this case ranged from 3 to 5. For the estimation we used a polynomial basis with $r=2$. 
To evaluate the performance, we computed the angle between the true reduction $\alphabf^T\Z$ and the estimated reduction $R(\X)$. We considered three choices for $R(\X)$: (i) the reduction is given by $\alphabfhat_{\tiny\hbox{PFC}}^T\X$, with $\alphabfhat_{\tiny\hbox{PFC}}$ computed as in standard PFC for continuous variables;
(ii) the reduction is $\alphabfhat_{\tiny\hbox{ORD}}^T\X$, with $\alphabfhat_{\tiny\hbox{ORD}}$ computed as proposed here but applied on the observed ordinal data $\X$;
(iii) $R(\X)=E(\alphabfhat_{\tiny\hbox{ORD}}^T\Z|\X)$ as proposed in Section 2.

Figure \ref{fig:simu}-\subref{bpnormal} shows boxplots of the obtained results for 100 runs of the experiment and a 
sample of size $n=500$.
The angle is measured in degrees.
It can be seen that the mean value of the angle is significantly smaller for the proposed method for ordinal predictors, compared to using choices (i) or (ii). 
It can be seen that the variance is somewhat increased, but the gain in accuracy clearly worths the price.  
It is less obvious from the figure that $\alphabfhat_{\tiny\hbox{ORD}}$ provides a better estimation of the true subspace spanned by $\alphabf$ than the standard $\alphabfhat_{\tiny\hbox{PFC}}$ estimator.
The 95\% normal confidence interval for the difference $[\hbox{angle}(\alphabf,\alphabfhat_{\tiny\hbox{PFC}})-\hbox{angle}(\alphabf,\alphabfhat_{\tiny\hbox{ORD}})]$ is $(9.59,11.02)$ degrees. These results show that for ordinal data, estimation of the subspace spanned by $\alphabf$ using the proposed method for ordered predictors clearly outperforms 
standard \texttt{PFC} as derived for continuous predictors.

It seems fair to ask whether the gain in performance discussed in this example still holds when the normality assumption for the underlying latent variables does not hold.
For standard  \texttt{PFC}, \citet[][Theorem 3.5]{cook_forzani_2008} showed that the estimator is still consistent when $\Z|Y$ deviates from multivariate normality.
To assess the performance of the proposed method in this scenario, we generated data similarly as before, but with $\epsilonbf$ non-normally distributed. 
In particular, we assumed that $\epsilonbf$ had chi-squared distributed coordinates $(\epsilonbf)_j\sim\chi^2(5)$. 
The rest of the simulation parameters remained fixed as before. 
The results obtained are shown in Figure~\ref{fig:simu}-\subref{bpchi2}. 
It can be seen that the angles obtained with both methods are very close to those obtained for conditionally normal data.  
This confirms the superiority of the proposed method and algorithm to estimate $\hbox{span}(\alphabf)$ when the predictors are ordered categories.

\begin{figure}
\begin{center}
   \subfigure[]{\includegraphics[width=0.5\textwidth]{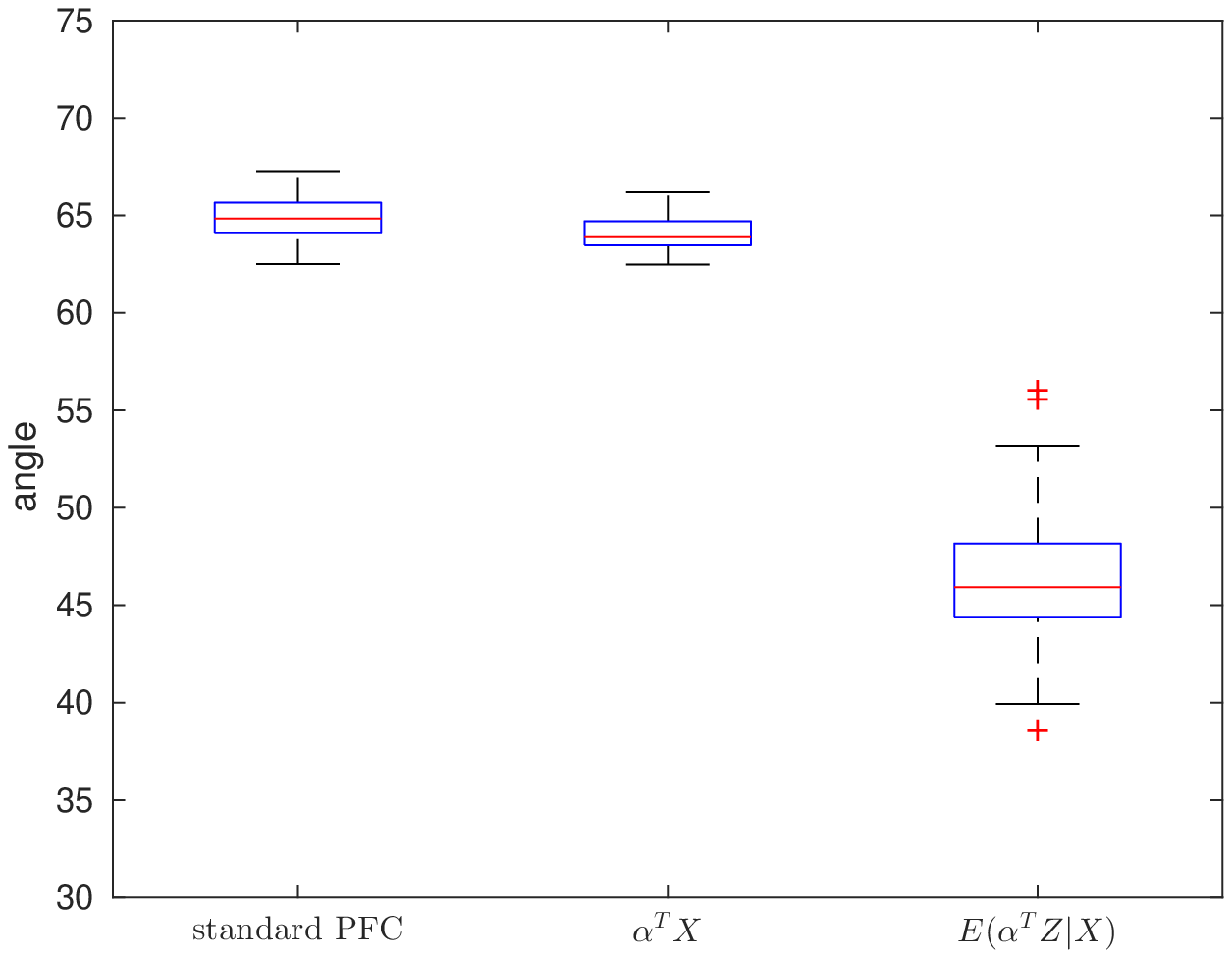}\label{bpnormal}}
   \hspace{-5mm}
   \subfigure[]{\includegraphics[width=0.5\textwidth]{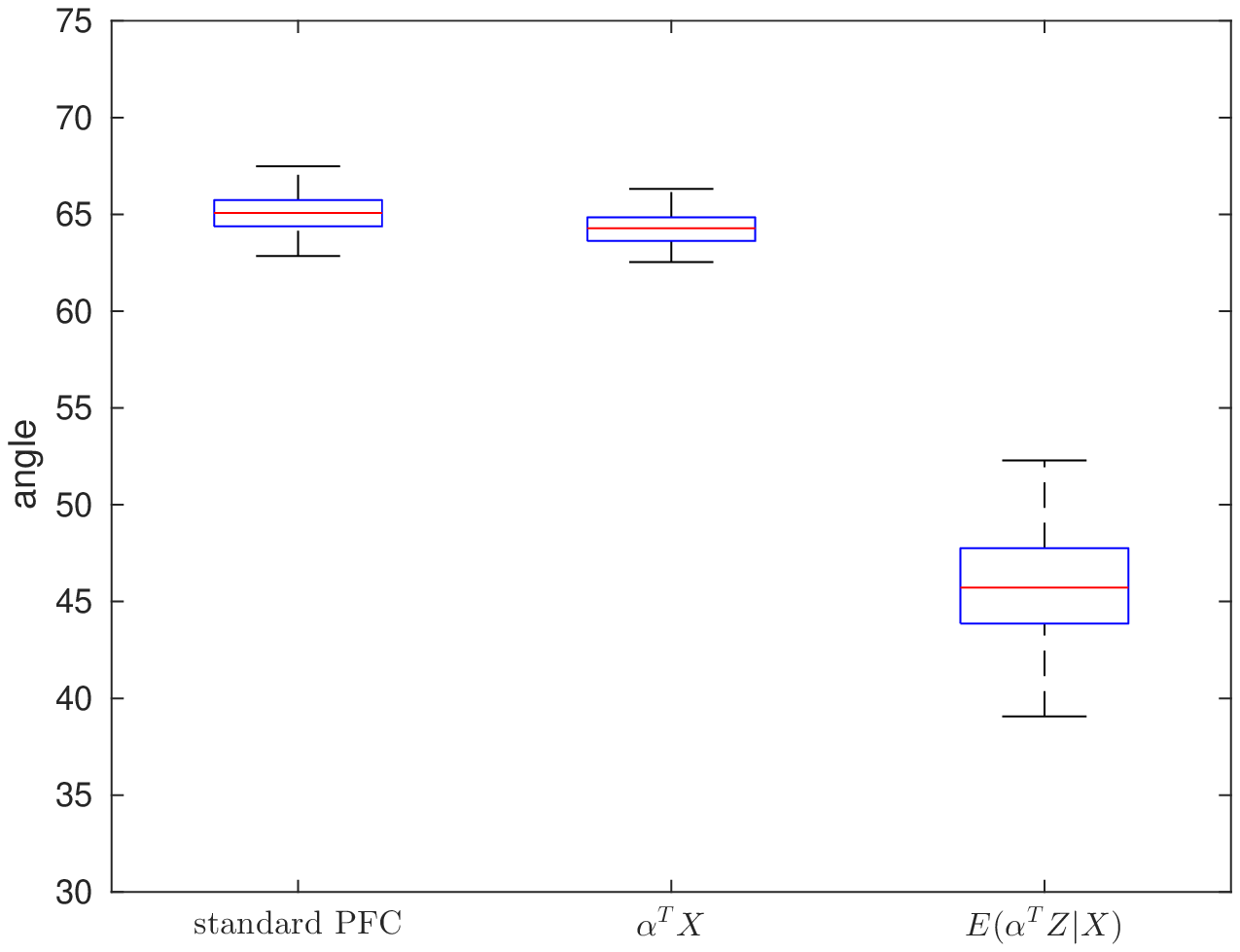}\label{bpchi2}}
\end{center}
\caption{Angle between $\hbox{span}(\alphabfhat)$ and $\hbox{span}(\alphabf)$, for $\alphabfhat$ 
    obtained using \texttt{PFC} and the proposed approach for ordinal predictors. 
\subref{bpnormal} normal conditional model for $\Z|(Y=y)$; \subref{bpchi2} non-normal conditional model for $\Z|(Y=y)$.}
\label{fig:simu}
\end{figure}

\subsection{Inference about $d$}
\textcolor{black}{
	In this section, we carry out a simulation study to evaluate methods to infer the dimension of the reduction subspace from the data. In particular, we compare the accuracy of permutation testing against 10-fold cross validation (CV) and well-known information criteria like Akaike's information criterion (AIC) and Bayes' information criterion (BIC).
	For this study, we generate data as described in Section \ref{sec:simu1}, with
	$p=10$, $d=2$ and we use a polynomial basis with degree $r=4$ for $\f_Y$.
	Since $d\leq\min(r,p)$, we search for the true value of $d$ within the
	set $\{0,1,2,3,4\}$. For permutation testing, we build the permutation distribution of
	$\hat{\Lambda}_m$ resampling the data 500 times as discussed in
	Section \ref{dimension} and use a significance level of $0.01$. For CV we used the averaged mean-squared prediction error over the test partition as the driving measure of performance, taking simple $k$-NN regression as the prediction rule.
	For information criteria, we take
	\begin{equation*}
	\hat{d}_{IC} = \mathop{\arg\min}_{d_o} \{2 \mathcal{Q}_{d_o}(\alphabf) + c_{IC}h_{\Theta}(d_o)\},
	\end{equation*}
	where $c_{AIC}=2$, $c_{BIC}=\log(n)$, $h_{\Theta}(d_o)=r d_o + d_o(p-d_o) + p(p+3)/2 + n_{\Theta}$ and $n_{\Theta}$ is the number of thresholds estimated during the computation of $\alphabfhat$.
	The experiment was repeated using two different sample sizes, $n=200$ and $n=300$, and it was run 500 times for each sample size. 
}

\textcolor{black}{
	Table \ref{tab:choose_d} shows the obtained results. 
	For $n=300$, permutation testing finds the true dimension 83\% of the runs, while CV finds it 66\% of the time. 
	In addition, for this sample size, both methods hardly ever pick less than two directions, and therefore no information is lost. The fraction of the runs that at most one extra direction is chosen, that is, $\hat{d}=2$ or $\hat{d}=3$, is 0.898 for permutation testing and 0.85 using cross-validation. On the other hand, information criteria show a poorer performance. AIC picks the right dimension around half of the runs. Nevertheless, the rest of the runs it picks $\hat{d}=1$, thus losing information. Underestimation of the true dimension is even more severe with BIC, since it picks only one direction almost always. 
}

\textcolor{black}{
	For the smaller sample size $n=200$, permutation testing is less
	accurate. It finds the true dimension 67\% of the runs in this
	scenario, but in 24\% of the runs it picked only one direction for projection instead of two. 
	CV finds the true dimension
	59\% of the runs, but tends to overestimate the required dimension
	$d$, leading to a potential loss in efficiency but preserving information.
	A test for the difference in the proportion of choices $\hat{d}\geq2$ between permutation testing and CV gives a $p$-value of $\approx 10^{-15}$, evidencing a statistical significant advantage of the cross-validation in this scenario of small samples in order to avoid information loss.
	On the other hand, in this setting AIC and BIC underestimate the true dimension even more frequently than for the larger sample size.
}

\textcolor{black}{
	Summarizing, both permutation testing and CV provide more accurate results compared to AIC and BIC, with information criteria typically underestimating the true dimension of the reduction. 
	Permutation testing seems to be a better procedure to infer $d$ when the sample size is large enough.
	Nevertheless, cross-validation can provide a safer solution regarding information loss when the available data is limited.}

\begin{table}[t]
	\centering
	\caption{Fraction of times a given value of $d$ was chosen}
	\begin{tabular}{lccccc}
		\hline
		&	     & \textsc{Permutation} & \textsc{CV} & AIC   & 	BIC \\
		\hline
		& $d=1$  & 0.240 & 0.000 &	0.794	&	0.998\\
		$n=200$ & $d=2$  & 0.670 & 0.591	&	0.206	&	0.002\\
		& $d=3$  & 0.067 & 0.214 	&	0	&	0\\
		& $d=4$  & 0.023 & 0.195 &	0 & 0\\
		\hline
		& $d=1$  & 0.077 & 0.000 & 0.488 & 0.986\\
		$n=300$ & $d=2$  & 0.832 & 0.657 & 0.512 & 0.014\\
		& $d=3$  & 0.066 & 0.191 & 0 & 0 \\
		& $d=4$  & 0.025 & 0.152 & 0 & 0\\
		\hline
		
	\end{tabular}
	\label{tab:choose_d}
\end{table}

\subsection{Performance of the proposed method including regularization}

Finally, we conducted a simulation study to assess the performance of the regularized version of the proposed method.
Unlike the previous setting, the reduction depends now on a subset of the predictors only. We chose that the first four predictors conveyed information
about the response, that is, 
$\alphabf = (\Abf \;\; \mathbf{0}_{2\times p-4})^T$, with
\begin{equation*}
\Abf = \left(
\begin{array}{cccc}
1/2 & 1/2 & 1/2 & 1/2\\
-1/2& 1/2 & 1/2 & -1/2
\end{array}\right).
\end{equation*}
We set the values for the rest of the parameters as described in Section \ref{sec:simu1}. 
The reduction was estimated using \texttt{PFC} for continuous predictors, the non-regularized method introduced in Section 4, and the regularized version proposed in Section 5.
In all cases we used a polynomial basis with degree $r=2$ for $\f_Y$.
For each estimator, we computed the angle between the subspaces $\mathcal{S}_{\alphabfs}$
and $\mathcal{S}_{\hat{\alphabfs}}$ in each of 100 runs of the
experiment. Figure \ref{fig:simu_lasso} shows the obtained results for $p=10$.
It can be seen that the obtained angles are smaller than in the case where all the predictors are relevant.
Since all the methods are applied to identical data, it is clear from the boxplot that the estimators specifically tailored to ordinal predictors still perform significantly better than the standard \texttt{PFC} approach. Moreover, for this situation, the regularized estimator (from now on, \texttt{reg-PFCord})  clearly proves to be superior to the non-regularized version. 

\begin{figure}[t!]
\begin{center}
    \includegraphics[width=0.65\textwidth]{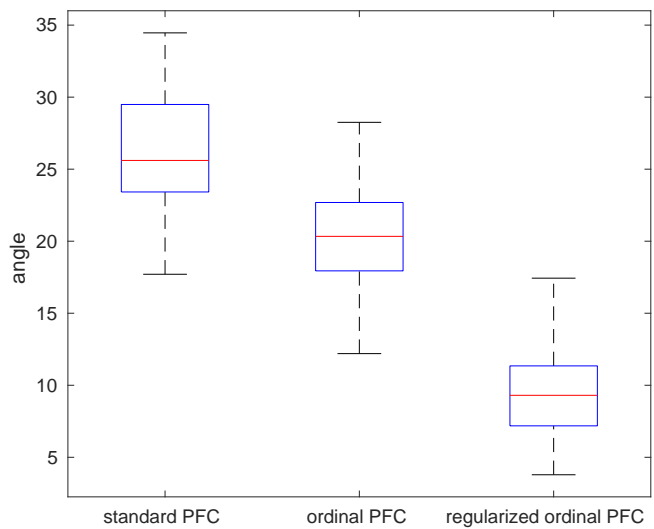}
\end{center}
\caption{Angle between $\hbox{span}(\alphabfhat)$ and $\hbox{span}(\alphabf)$, for $\alphabfhat$ obtained using \texttt{PFC}, the proposed approach for ordinal predictors without regularization, and the regularized approach introduced in Section 5.}
\label{fig:simu_lasso}
\end{figure}

To further study the performance of the proposed regularized estimator, another set of experiments was carried out in order to evaluate the stability of the subset of variables chosen by the algorithm. Denote by $S_0 \subset \{1,2,\dots,p\}$ the index set for the subset of variables that are truly relevant for describing the response $Y$ and let $S_{0}^c$ be its complement. Similarly, let $\Shat$ be the subset of variables chosen by the regularized
estimator, in the sense that $\alphabf_j=\mathbf{0}$ for $j \notin \Shat$, and let $\Z_{S_0}$ be the random
vector with entries $Z_j$ for $j\in S_0$, with a similar definition for $\Z_{S^c_0}$.
We are interested in assessing: (i) $\Pr(S_0 \subset \Shat)$, as an indicator
of the relevant variables that are indeed retained; (ii) the average cardinality of the
set $\Shat$ ($\#\Shat$) as a measure of the amount of non-relevant variables
that are preserved. Moreover, we are interested in evaluating how these performance measures vary for different amounts 
of correlation between $\Z_{S_0}$ and $\Z_{S^c_0}$. To measure this dependence we use the distance correlation measure as defined in \citep{dcov2007,dcov2009}, denoted here by 
$dCor_n(\Z_{S_0},\Z_{S^c_0})$ when it is computed from a sample of size $n$. This is a generalized nonparametric measure of correlation that is suitable for random vectors of different size and it does not require tuning any parameter for its computation.
To control this quantity, we adjust the value of $\Deltabf$ used to generate the data.
In particular, we set $\Deltabf = 4\I_p + \rho\alphabf\B\alphabf^T$ for
different values of $\rho$. We obtained $dCor_n(\Z_{S_0},\Z_{S^c_0})\approx 0.3$ for $\rho=0.2$, $0.5$ for $\rho=0.3$ and approximately $0.7$ for $\rho=0.5$. During all the experiments, the number of relevant variables in $S_0$ was held fixed at 4, with $p=20$. 

\begin{table}[t!]
\centering
\caption{Performance of variable selection algorithms when using \texttt{reg-PFCord} or \texttt{reg-PFC}.}

\vspace{3mm}

\begin{tabular}{c c c c c c}
    \hline
    $n$  	& $\rho$ & \multicolumn{2}{c}{\texttt{reg-PFCord}}& \multicolumn{2}{c}{\texttt{reg-PFC}} \\
    &  & $\Pr(S_0 \subseteq \Shat)$ & $\#\Shat$ & $\Pr(S_0 \subseteq \Shat)$ & $\#\Shat$\\
    \hline
     & 0.0	& 0.98 & 3.87 ($\pm 0.15$)&	0.92& 3.88 ($\pm 0.63$)\\
    200	& 0.2   & 0.90 & 3.91 ($\pm 0.29$)&  0.13& 2.36 ($\pm 0.76$)\\
     & 0.3	& 0.53	& 3.49 ($\pm 0.61$)&  0.04& 2.23 ($\pm 0.46$)\\
     & 0.5   & 0.21 & 2.86 ($\pm 0.63$)&  0.01& 2.07 ($\pm 0.38$)\\
    \hline
     & 0.0  & 1.00 & 4.11 ($\pm 0.00$)  &  1.00 & 4.00 ($\pm 0.00$)\\
    500	& 0.2   & 0.98& 4.05 ($\pm 0.10$)  & 0.96 & 4.52 ($\pm 0.73$)\\
     & 0.3   & 0.96 & 3.92 ($\pm 0.38$)  & 0.63 & 3.52 ($\pm 0.98$)\\
       	& 0.5   & 0.54 & 4.05 ($\pm 1.08$)  & 0.17 & 2.53 ($\pm 0.83$)\\
    \hline
\end{tabular}
\label{tab:vs2}
\end{table}

Table \ref{tab:vs2} shows the results obtained using 100 replicates of the
experiment for each assessed condition.
It can be seen that for a large enough sample ($n=500$), the penalized versions of both
\texttt{PFCord} and \texttt{PFC} achieve perfect accuracy in selecting the true
active set when the predictors are not correlated each other.
For moderate levels of correlation between the predictors ($\rho=0.3$), using \texttt{reg-PFCord} the true active set $S_0$ is contained in the solution 96\% of the
time, with a very low fraction of false insertions. 
On the other hand, using the regularized version of the standard \texttt{PFC} (from now on, \texttt{reg-PFC}) in the same scenario allows picking the true active set of
variables only 63\% of the time. 
This difference is statistically significant at the 0.001 level. 
Moreover, the average number of variables picked by the \texttt{reg-PFC} is around 3.50, meaning that some information is
typically lost in the procedure.
For very high levels of correlation between the predictors ($\rho=0.5$),
the true set of active variables is picked only 17\% of the time using
\texttt{reg-PFC}, whereas \texttt{reg-PFCord} still finds it half the time. 
The loss of accuracy of \texttt{reg-PFCord} usually involved replacing one of the predictors in the true active set by a highly correlated alternative, mantaining the average cardinality of the estimated set $\Shat$
close to 4.
When the sample size is smaller ($n=200$), the performance of \texttt{reg-PFC} for variable selection degrades much faster with the
level of correlation than the ordinal counterpart. 
When the predictors are
uncorrelated, \texttt{reg-PFC} picks the true active set $S_0$ 92\% of the
time, whereas \texttt{reg-PFCord} does it 98\% of the time. 
But for low levels
of correlation between the predictors ($\rho=0.2$), the performance of {\texttt{reg-PFC}} decreases quickly to 13\%, while for the ordinal counterpart it is
still greater than 90\%. 
For very high
levels of correlation, both procedures tend to underestimate the number of
relevant variables, with this trend being stronger for the { \texttt{reg-PFC}.}
These results show that using the proposed method especially tailored to ordinal
data provides significantly higher accuracy when variable selection is needed.

\section{Real data analysis: SES index construction} \label{realdata}
%
%

For many social protection and welfare programs carried out by governments and NGOs, a 
cla\-ssification of households or individuals into different socio-economic groups is required. 
For example, over the last decades, it has become common in developing countries that governments establish economic aid programs focused on the most deprived households. 
This tailoring of the aid is achieved via a focalization index, which basically mounts to a Socio-Economic Status ({\em SES}) index, as it is commonly known in the related literature. 
In particular, in many Latin American countries, this focalization index (called {\em \'Indice de Focalizaci\'on de Pobreza}) has been used to implement several programs to reduce poverty (e.g., the CAS in Chile, Sisben in Colombia, SISFOH in Per\'u, Tekopor\'a in Paraguay, SIERP in Honduras, and PANES in Uruguay, among others).

Income or consumption expenditures constitute a traditional focus of poverty analysis, and some 
countries take an income-based poverty line from a household survey to infer the socioeconomic 
situation of the population \citep{Mokamane12,Richardson12}. 
However, the collection 
of income data presents many problems in terms of unavailability or unreliability \citep{Vyas06,Doocy06}. For this reason, asset-based indexes are often constructed as a proxy of income, taking into account some housing and household variables that are easier to observe. 
This proxy, usually called the SES index, is most of the time computed using principal component analysis (\texttt{PCA}) \citep{Merola14,Hoque14}. 
Since the observable variables used to construct these indexes are ordered categorical  variables, \cite{Kolenikov09} proposed a variant of \texttt{PCA} adapted for ordinal data using polychoric correlations between predictors instead of the standard covariance matrix.

In this example, we provide a different approach to the construction of an SES index, based on the proposed SDR methodology for ordinal data. 
The main idea is to obtain a single index to predict a unidimensional measure of some socioeconomic aspect, such as household income or the poverty condition. Therefore, in this case we fix the dimension of the reduction to be $d=1$ and derive the index as a normalized version of the \textcolor{black}{supervised dimension reduction $\alphabfhat^T\E(\Z|\X) \in \real$}. 
Unlike \texttt{PCA}-based indexes, this new approach uses information about the response under analysis.

The data comes from the microdata of the \textit{Encuesta Permanente de Hogares}  (EPH) of Argentina, 
taking the fourth trimester of 2013. 
The EPH is the main household survey in Argentina and is 
carried out by the \textit{Instituto Nacional de Estad\'isticas y Censos} (INDEC). 
We consider nine 
ordinal variables about household living conditions, and two socio-economic variables of heads of 
households (educational attainment and work situation). 
More details about these variables can be 
found in  \ref{appendiceD}. 
To take into account regional heterogeneity, we estimate separate SES 
indexes for the following five regions: the metropolitan area of Buenos Aires ($n=2351$ 
households), Humid Pampas ($n=5003$), the Argentine Northwest ($n=2852$), the Northeast 
($n=1594$), and Patagonia ($n=2398$).
Two cases with different types of response are considered: a continuous one, household income per capita (\texttt{ipcf}), and 
a binary one based on income (\texttt{poverty}) that indicates whether a household is poor or not. We are interested in demonstrating that the proposed \texttt{reg-PFCord} provides a superior alternative to a \texttt{PCA}-based method for constructing an SES index, while retaining a predictive power comparable to the full set of predictors. 
To do this, the predictive performance of the proposed response-driven index is compared to the following strategies:
\begin{itemize}
\item The full set of predictors is included without dimension reduction and they are treated as continuous (metric) predictors. 
We will refer to this apporach as \texttt{FULL}.
\item The full set of predictors is included without dimension reduction and they are treated through dummy variables. 
We will refer to this approach as \texttt{FULL-I}.
\item The full set of predictors is considered but using a group-lasso-type procedure for ordinal predictors that induces variable selection \citep{Gertheiss10}. 
We will refer to this approach as \texttt{LASSOord}.

\item A variant of \texttt{PCA} tailored to ordered categorical predictors using polychoric correlations \citep{Kolenikov09}. 
We will refer to this method as \texttt{PCApoly}.
\item A nonlinear variant of \texttt{PCA} that uses special scaling to take into account the ordered categories \citep{Linting09}. 
We will refer to this approach as \texttt{NLPCA}.
\item Standard moments-based sufficient dimension reduction methods \texttt{SIR}, \texttt{SAVE} and \texttt{DR}.
\end{itemize}

The first three strategies are included in order to provide a reference for the performance achievable using the full set of predictors, but it should be clear that they do not provide an index. 
Actually, only the last three approaches in the list are competing methods for SES-index extraction. 
Among them, \texttt{PCApoly} and \texttt{NLPCA} can deal explicitely with ordinal predictors and they will be further compared later. 

For every strategy, we fit a logistic regression for the \texttt{poverty} response and a linear regression for the  \texttt{ipcf} response. 
When computing the reduction, we use a different choice of $\f_Y$ for each type of outcome. For the continuous response, we use a polynomial basis with degree $r=2$. 
For the binary response, $\f_Y$ is simply a centered indicator variable. The data was partitioned into ten disjoint sets to allow for ten replications of the experiment. In each run, one of the subsets was used as the test set, while the rest of them formed the training sample. Averaged 10-fold cross-validation MSEs obtained with each method are shown in Table \ref{tab:eph}, along with the corresponding standard deviations.

\begin{table}[t]
\small{
    \centering
    \caption{10-fold MSE for SES index}

    \begin{tabular}{rrrrrrr}
        \hline

        &                                 \multicolumn{ 6}{c}{Prediction Errors -MSE} \\
        
        Response  &     Method & {\it Buenos Aires} & {\it Humid Pampas} & {\it Northwest} & {\it Northeast} & {\it Patagonia} \\
        \hline
         {\it Per capita Income} 
        
         & \textsc{\texttt{reg-PFCord}} & { 7.29} & { 4.72} & { 4.73} & {3.34} & { 12.80} \\
       (continuous)    &         & { (2.91)} & { (1.73)} & { (2.69)} & { (1.52)} & {  (3.72)} \\            
&    \texttt{PCApoly} &      7.60 &      5.10 &      5.07 &      3.68 &      14.7 \\
        &            &    (2.45) &    (0.90) &    (1.77) &    (0.90) &    (4.01) \\
        & \texttt{NLPCA} & 7.38 & 4.95 & 4.89& 3.52 & 13.67\\
        &            &    (2.29) &    (0.61) &    (1.48) &    (0.65) &    (3.71) \\
        \cline{2-7}
        &\texttt{SIR}	&	7.36	&	6.21	&	5.61	&6.15	&	14.41	\\
        &				&	(4.38)	&	(4.73)	& (3.72)	&	(0.86)	&	(1.16)\\
        &\texttt{SAVE}	&  9.06	&	6.09	&	5.87	&	4.12	&	16.21 \\
        &				& (4.04)	& (0.99)	&	(2.73)	&	(1.24)	&	(3.74)	\\
        &\texttt{DR}		& 8.96	&	5.76	&	5.85	&	4.10	&	15.95	\\
        &				&	(3.96)	&	(1.02)	&	(2.74)	&	(1.24)	&	(3.70)	\\
        \cline{2-7}        
        & \texttt{FULL} &      7.22 &      4.69 &      4.68 &      3.32 &      13.14 \\
 &            &    (3.50) &    (0.88) &    (2.47) &    (0.93) &    (3.34) \\
        & \texttt{FULL-I} & 7.01 & 4.52 & 4.48 & 3.08 & 12.92\\
        &            &    (2.46) &    (0.83) &    (1.74) &    (0.76) &    (3.80) \\
        & \texttt{LASSOord} & 7.00& 4.51 & 4.42& 3.05& 12.88 \\
        &                  &    (2.46) &    (0.83) &    (1.77) &    (0.76) &    (3.84) \\

        &            &            &            &            &            &            \\
         \hline
          
        {\it Poverty} &  \texttt{reg-PFCord} 
        & { 0.204} & { 0.169} & { 0.278} & { 0.288} & { 0.126} \\
       (discrete)   &          & { (0.017)} & { (0.012)} & { (0.031)} & { (0.029)} & { (0.025)} \\

         &    \texttt{PCApoly} &     0.213 &     0.188 &     0.324 &     0.357 &     0.133 \\
        &            &    (0.024) &    (0.020) &    (0.026) &    (0.053) &    (0.027) \\
        & \texttt{NLPCA} & 0.212 & 0.188 & 0.325& 0.358 & 0.134 \\
        &            &    (0.023) &    (0.019) &    (0.025) &    (0.055) &    (0.027) \\        \cline{2-7}
        &\texttt{SIR} & 0.229	& 0.204	&	0.366	&	0.392	& 0.1314	\\
				      && (0.018) & (0.010) & (0.032) & (0.047) & (0.022) \\
		&\texttt{SAVE}& 0.229	&	0.204	&	0.362	&	0.378	&	0.133	\\
		&			  & (0.019)	&	(0.009)	&	(0.031)	&	(0.042)	&	(0.022)\\
		&\texttt{DR}  & 0.230	& 	0199	&	0.357	&	0.364	&	0.133	\\
		&			&	(0.019)	&	(0.009)	&	(0.035)	&	(0.040) & 	(0.023)\\	            
        \cline{2-7}
        &\texttt{FULL} &     0.202 &     0.162 &     0.274 &     0.287 &     0.129 \\
        &            &    (0.021) &    (0.008) &    (0.026) &    (0.036) &    (0.020) \\
        & \texttt{FULL-I} & 0.206 & 0.171 & 0.286& 0.298 & 0.126 \\
        &            &    (0.023) &    (0.021) &    (0.024) &    (0.065) &    (0.028) \\
         & \texttt{LASSOord} & 0.206 & 0.171 & 0.286 & 0.298 & 0.129\\
        &            &    (0.023) &    (0.021) &    (0.024) &    (0.065) &    (0.027) \\
        &            &            &            &            &            &            \\
        \hline
        \multicolumn{5}{c}{ Note: standard deviations in parentheses. Database: EPH (2013)}\\
    \end{tabular}  
    
    \label{tab:eph}}
\end{table}

From the table it can be seen that for the continuous response, using dummy variables for the full set of predictors, as in \texttt{FULL-I} and \texttt{LASSOord}, is more effective than considering the full set of predictors as continuous variables (\texttt{FULL}).
Among the SES indexes, scores show that \texttt{reg-PFCord} is superior to both \texttt{PCApoly} and \texttt{NLPCA}, with the latter being slightly superior to the former. In addition, all these methods specifically targeted to ordinal predictors perform better than standard sufficient dimension reduction methods represented here by \texttt{SIR}, \texttt{SAVE} and \texttt{DR}, which were originally aimed to continuous predictors only. Moreover, prediction errors obtained with \texttt{reg-PFCord} are very close to those attained with \texttt{FULL} across all the regions. 

On the other hand, it is interesting to see that for the binary outcome, unlike the continuous case, \texttt{FULL} performs better than \texttt{FULL-I} and \texttt{LASSOord}.
Among the SES indexes, the predictive performance of \texttt{reg-PFCord} is again very close to or identical with that of \texttt{FULL}, and it outperforms the PCA-based methods and standard moment-based sufficient dimension reduction methods in this scenario too. 
Moreover, \texttt{reg-PFCord} outperforms \texttt{LASSOord} for three of the regions when the discrete response is considered. It should also be remarked that, unlike \texttt{LASSOord}, obtaining the indexes from the \textcolor{black}{SDR-based techniques} 
allows us to use any predictive method.

As an illustration of the obtained fit, Figure \ref{fig:sesplots} shows marginal model plots for the regression of \texttt{ipcf}  on the SES index obtained for the whole database. A quadratic term $\hbox{SES}^2$ was added 
to correct for curvature in the estimated regression function and the response was transformed by 
$\hbox{ipcf}\leftarrow\hbox{ipcf}^{1/3}$ following a Box--Cox transformation analysis. It can be seen that for SES modeled using \texttt{PCA} on 
polychoric correlations, the index values are concentrated mainly in the interval $[0.5;1.0]$ 
whereas, for SES modeled using \texttt{reg-PFCord}, the spread of the index values is more regular over the whole 
interval $[0;1.0]$. 
This allows for a better fit of the linear model, as shown by an $R^2$ value of 
$0.302$ compared to $0.231$ obtained with the SES index based on \texttt{PCA}.

\begin{figure}
\centering
\includegraphics[width=.45\textwidth,height=2in]{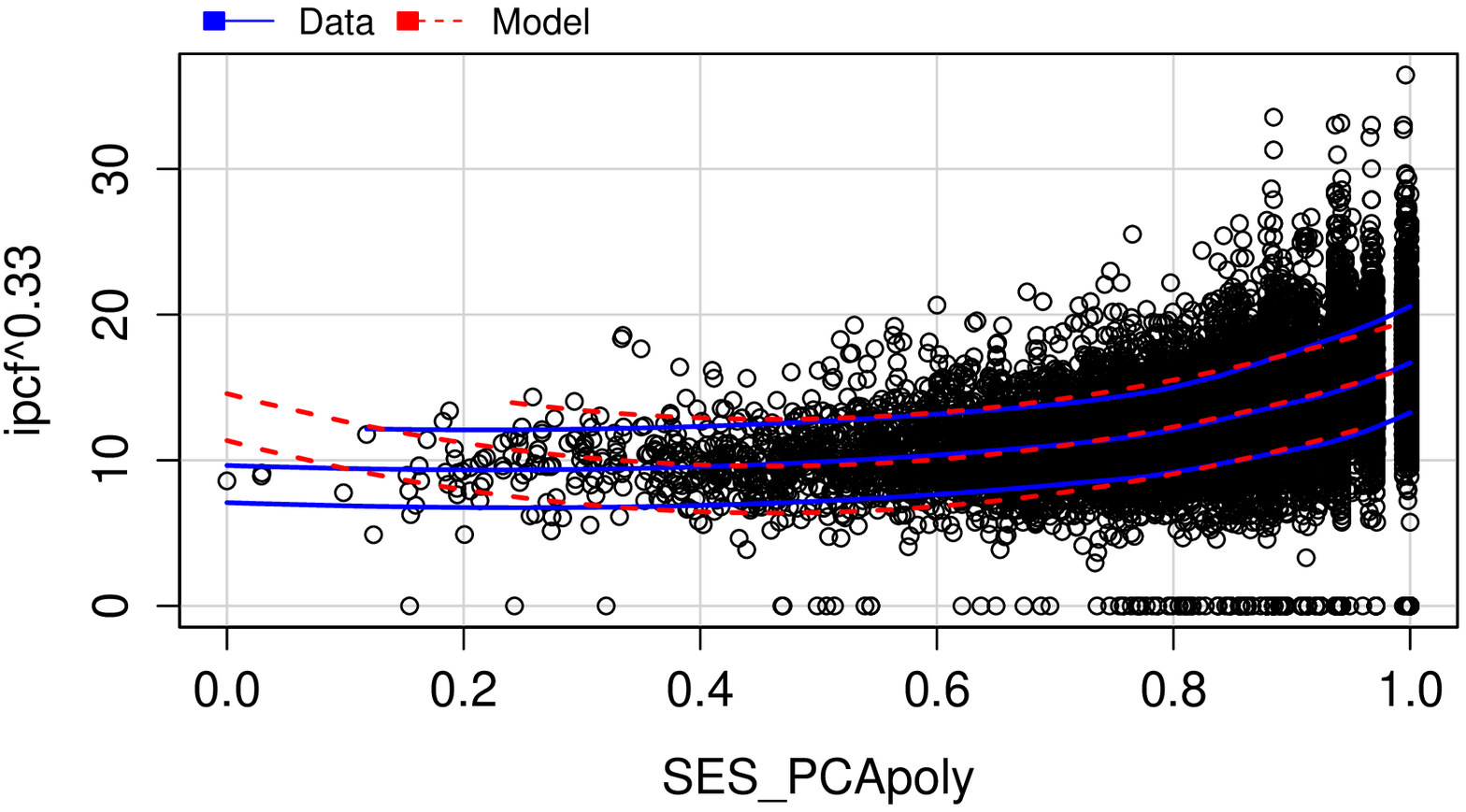}
\includegraphics[width=.45\textwidth,height=2in]{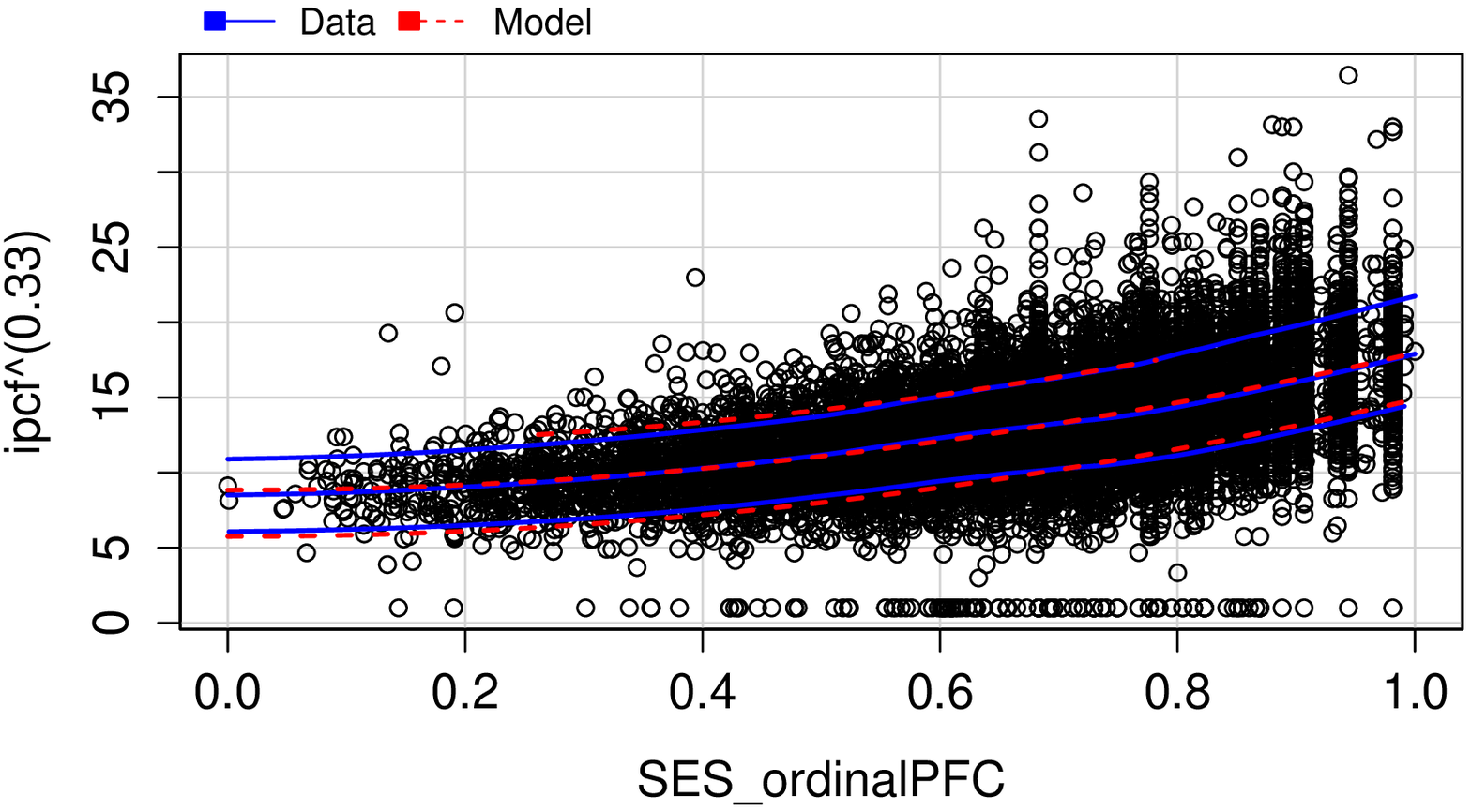}
\caption{Marginal model plots showing the fit of linear model of income as a function of the 
obtained SES index. 
\texttt{PCA}-based SES index is shown on the left and the 
proposed method on the right.}\label{fig:sesplots}
\end{figure}

Tables \ref{tab:coef} and \ref{tab:coef2} show the estimated coefficient vectors that 
define the SES index using \texttt{ipcf} and \texttt{poverty} as response variables, respectively. 
To keep the analysis clear, standard methods not targeted to ordinal predictors like \texttt{SIR}, \texttt{SAVE} and \texttt{DR} were not included, since they showed a clearly inferior performance for prediction in Table \ref{tab:eph}.
Note that for the proposed method, some of the elements of $\alphabfhat$ have been pushed to zero in the regularized estimation, whereas for \texttt{PCApoly} and \texttt{NLPCA} only \textit{working hours} seems not to be relevant for constructing the index. In addition, several differences can be appreciated between the \texttt{reg-PFCord} and \texttt{PCA}-based approaches (i.e., \texttt{PCApoly} and \texttt{NLPCA}) from the reported results. First, the relative importance of each predictor in the 
SES index obtained is different for the two methods. For instance, \textit{overcrowding} obtains the highest weight with \texttt{reg-PFCord} across all the regions for both responses, whereas \textit{toilet facility} and \textit{water location} appear as the most important in index construction based on \texttt{PCApoly}, and \textit{toilet facility} and \textit{toilet drainage} for \texttt{NLPCA}. Second, SES indexes constructed using both \texttt{PCA} methods give similar weights to the predictors across the different regions. On the other hand,  SES indexes based on \texttt{reg-PFCord} capture the regional economic divergence explained by different factor endowments, productivity, 
activity levels and regional economic growth patterns. Moreover, in the richest 
Argentinian urban regions (specifically, Buenos Aires and Humid Pampas) the regularized estimation 
of \texttt{reg-PFCord} often sets to zero the variables with more weight in \texttt{PCA}-derived SES index. This difference is appealing, since these regions have in general  better services and public infrastructure, so that variables related to drainage, source of water and toilet facility are less important for measuring socio-economic status. 
In this way, other variables, such as overcrowding or schooling, are needed in order to have a better SES index to predict household income. In the same line, for regions with higher levels of poverty (Northwest and Northeast) the \texttt{reg-PFCord}-based SES index shows that other variables, such as \textit{housing location}, \textit{source of drinking water} or 
\textit{water location}, become important for determining socio-economic status.

Comparing both \texttt{PCA} methods, it can be noted that the \texttt{NLPCA} is more sensitive to the regional heterogeneity than is the \texttt{PCApoly}, but differences in the index weights compared to those of the \texttt{reg-PFCord} remain substantial. Additionally, it can be appreciated that SES indexes obtained using \texttt{reg-PFCord} are sensitive to the response variable used to characterize a social phenomenon of interest. 
For example, in Buenos Aires, Humid Pampas and Patagonia, \textit{schooling} has a considerable weight in the SES index to explain per capita income but not to predict poverty. 
This makes sense, since for these richest regions it is easier for all the population to get access to basic levels of schooling. On the other hand, the decision to pursue higher levels of education is often driven by income. 
Moreover, for these richest regions, some 
variables, such as \textit{toilet drainage}, \textit{toilet facility} or \textit{toilet sharing} become 
relevant to explaining whether a household is poor or not (following poverty line criteria). Such differences cannot be captured by an SES index based on \texttt{PCA} methods.

\begin{landscape} 
\begin{table}[t!]
\centering
\caption{Comparison of SES index results for \texttt{ordinal PCF} and \texttt{PCApoly} to predict household per capita income.}
\footnotesize
\begin{tabular}{r|rrr|rrr|rrr}
\hline
Variables & \multicolumn{ 3}{|c}{Buenos Aires} & \multicolumn{ 3}{|c}{ Humid Pampas} & \multicolumn{ 3}{|c}{Northwest} \\
\hline
       &\texttt{reg-PFCord} &  \texttt{PCApoly} & \texttt{NLPCA}  
       & \texttt{reg-PFCord} &   \texttt{PCApoly} & \texttt{NLPCA}  
       &\texttt{reg-PFCord} &   \texttt{PCApoly} &   \texttt{NLPCA} \\
\hline
{\it housing location} 
&          0 &    -0.1690 & -0.0943     
&    0 		&    -0.1903 &    -0.0976
&-0.1314   &    -0.1068 &    -0.0835\\

{\it housing quality} 
&    -0.1985 &    -0.3768 & -0.2199
&     0.2591 &    -0.3557 & -0.1985   
&    0 &    -0.3278 & -0.1849\\

{\it sources of cooking fuel} 
&    -0.4646 &    -0.3788 &  -0.2080   
&     0.3627 &    -0.3609 & -0.1678 
&    -0.1070 &    -0.3287 & -0.1582\\

{\it overcrowding} 
&    -0.7272 &    -0.2888 & -0.1788
&     0.8300 &    -0.2351 & -0.1329  
&    -0.8798 &    -0.1991 & -0.1194\\

{\it schooling} 
&    -0.2676 &    -0.2275 & -0.1474           
&     0.2668 &    -0.2075 & -0.1135
&    -0.3614 &    -0.2197 & -0.1201\\

{\it toilet drainage} 
&    -0.0873 &    -0.3381 &  -0.2047 
&          0 &    -0.3519 &  -0.2333  
&    -0.0901 &    -0.3623 & -0.2462\\

{\it toilet facility} 
&          0 &    -0.4061 & -0.2246 
&     0 &    -0.4105 & -0.2411
&     0 &    -0.4217 & -0.2545\\

{\it toilet sharing} 
&          0 &    -0.2759 & -0.1186
&          0 &    -0.3176 & -0.1699
&          0 &    -0.2579 & -0.1334\\

{\it water location} 
&     0.1700 &    -0.3918 & -0.1790 
&          0 &    -0.3933 & -0.1941 
&    -0.2054 &    -0.4202 & -0.2309\\

{\it water source} 
&          0 &    -0.2023 & -0.1033
&          0 &    -0.2461 & -0.1129
&    0 &    -0.3646 & -0.1374\\

{\it working hours} 
&    -0.3283 &     0      & 0
&     0.2029 &     0      & 0
&    -0.1277 &     0      & 0\\

\hline
& \multicolumn{ 3}{|c}{Northeast} & \multicolumn{ 3}{|c|}{Patagonia} & \\
\hline 
&\texttt{reg-PFCord} &  \texttt{PCApoly} &  \texttt{NLPCA} 
       &\texttt{reg-PFCord} & \texttt{PCApoly} & \texttt{NLPCA} & \\
\hline           
{\it housing location} 
&-0.1509 &    -0.1809 &    -0.0978
&-0.1149 &    -0.1437 &  -0.0981\\

{\it housing quality} 
&    0 &    -0.3727 & -0.2130
&    -0.3046 &    -0.3258 & -0.1844\\

{\it sources of cooking fuel} 
&    -0.0742 &    -0.1648 & -0.0646
&    -0.1797 &    -0.4026 &  -0.1810\\

{\it overcrowding} 
&    -0.8496 &    -0.2052 & -0.1040
&    -0.7263 &    -0.2207 & -0.0984\\

{\it schooling} 
&    -0.3507 &    -0.1869 & -0.1009
&    -0.3670 &    -0.1284 & -0.0516\\

{\it toilet drainage} 
&          0 &    -0.3572 &  -0.2573
&    -0.1383 &    -0.4122 & -0.2566\\

{\it toilet facility} 
&    0 &    -0.4383 & -0.2735
&          0 &    -0.4376 & -0.2622\\

{\it toilet sharing} 
&    -0.2284 &    -0.2921 & -0.1344
&    -0.1204 &    -0.2937 & -0.1734\\

{\it water location} 
&    0 &    -0.4227 & -0.2377
&     0.1877      &    -0.4169 &  -0.2196\\

{\it water source} 
&    -0.2574 &    -0.3733 & -0.1384
&     0.2473 &    -0.1525 & -0.0522\\

{\it working hours} 

&    -0.0922 &     0 &  0
&    -0.2637 &     0 & 0\\
\hline
\label{tab:coef}
\end{tabular}  
\end{table}
\end{landscape}
\begin{landscape} 
\begin{table}[t!]
\centering
\caption{Comparison of SES index results for \texttt{ordinal PCF} and \texttt{PCApoly} for a discrete response (poverty).}
\footnotesize

\begin{tabular}{r|rrr|rrr|rrr}
\hline
Variables & \multicolumn{ 3}{|c}{Buenos Aires} & \multicolumn{ 3}{|c}{ Humid Pampas} & \multicolumn{ 3}{|c}{Northwest} \\
\hline
       &\texttt{reg-PFCord} &  \texttt{PCApoly} & \texttt{NLPCA}  & \texttt{reg-PFCord} &   \texttt{PCApoly} & \texttt{NLPCA}  &\texttt{reg-PFCord} &   \texttt{PCApoly} &   \texttt{NLPCA} \\
\hline

{\it housing location} 
&          0 &    -0.1690 & -0.0943
&          0 &    -0.1903 & -0.0976
&    -0.2434 &    -0.1068 & -0.0835\\

{\it housing quality} 
&    -0.4033 &    -0.3768 & -0.2199
&     0.3347 &    -0.3557 & -0.1985
&    0       &    -0.3278 & -0.1849\\

{\it sources of cooking fuel} 
&    -0.5240 &    -0.3788 & -0.2080
&     0.3579 &    -0.3609 & -0.1678  
&    0       &    -0.3287 & -0.1582\\

{\it overcrowding} 
&    -0.7076 &    -0.2888 & -0.1788
&     0.7216 &    -0.2351 & -0.1329
&    -0.7939 &    -0.1991 & -0.1194\\

{\it schooling} 
&          0 &    -0.2275 & -0.1474
&          0 &    -0.2075 & -0.1135
&    -0.2094 &    -0.2197 & -0.1201\\

{\it toilet drainage} 
&    0       &    -0.3381 & -0.2047
&          0 &    -0.3519 & -0.2333
&          0 &    -0.3623 & -0.2462\\

{\it toilet facility} 
&    0 &    -0.4061 & -0.2246
&     0.3990 &    -0.4105 & -0.2411
&    -0.1528 &    -0.4217 & -0.2545\\

{\it toilet sharing} 
&    -0.1836 &    -0.2759 & -0.1186
&     0		 &    -0.3176 & -0.1699
&          0 &    -0.2579 & -0.1334\\

{\it water location} 
&    -0.1208 &    -0.3918 & -0.1790
&     0.2647 &    -0.3933 & -0.1941
&    -0.4933 &    -0.4202 & -0.2309\\

{\it water source} 
&          0 &    -0.2023 & -0.1033
&          0 &    -0.2461 & -0.1129
&          0 &    -0.3646 & -0.1374\\

{\it working hours} 

&    -0.1173 &     0	 & 0
&     0.0990 &     0	 & 0
&          0 &     0	 & 0\\
\hline
& \multicolumn{ 3}{|c}{Northeast} & \multicolumn{ 3}{|c|}{Patagonia} & \\
\hline 
&\texttt{reg-PFCord} &  \texttt{PCApoly} &  \texttt{NLPCA} 
       &\texttt{reg-PFCord} & \texttt{PCApoly} & \texttt{NLPCA} & \\
\hline        
{\it housing location} 
&    -0.1982 &    -0.1809 & -0.0978
&    -0.1187 &    -0.1437 & -0.0981\\

{\it housing quality} 
&          0 &    -0.3727 & -0.2130
&    -0.3693 &    -0.3258 & -0.1844\\

{\it sources of cooking fuel} 
&    -0.2509 &    -0.1648 & -0.0646
&    -0.2788 &    -0.4026 & -0.1810\\

{\it overcrowding} 
&    -0.7063 &    -0.2052 & -0.1040
&    -0.3987 &    -0.2207 & -0.0984\\

{\it schooling} 

&    -0.1442 &    -0.1869 & -0.1009
&    -0.0766 &    -0.1284 & -0.0516\\

{\it toilet drainage} 
&          0 &    -0.3572 & -0.2573
&    -0.1887 &    -0.4122 & -0.2566\\

{\it toilet facility} 
&    0		 &    -0.4383 & -0.2735
&    -0.1313 &    -0.4376 & -0.2622\\

{\it toilet sharing} 
&    -0.3477 &    -0.2921 & -0.1344
&    -0.1289 &    -0.2937 & -0.1734\\

{\it water location} 
&     0      &    -0.4227 & -0.2377
&     0.2785 &    -0.4169 & -0.2196\\

{\it water source} 
&    -0.5071 &    -0.3733 & -0.1384
&     0.6585 &    -0.1525 & -0.0522\\

{\it working hours} 

&          0 &     0	 & 0
&    -0.1626 &     0	 & 0\\
\hline
\label{tab:coef2}
\end{tabular}  
\end{table}
\end{landscape}

\section{Conclusions}
\label{Conclusion}

The approximate expectation-maximization (EM) algorithm presented here for dimension reduction in 
regression problems with ordinal predictors is proved to outperform the standard inverse regression methods derived for continuous predictors, both in simulation settings and with real data sets involving ordered categorical predictors.
Experiments showed that this advantage is emphasized in variable selection applications, where the proposed method clearly outperforms its counterpart for continuous data when the counterpart is naively applied to ordinal predictors.
This is not a minor issue since many analyses in the applied sciences usually treat them as continuous variables, not taking into account their discrete nature. Moreover, it has better computing efficiency due to the proposed approximate EM algorithm's rendering the method feasible for a much larger set of problems compared to using the exact computation of the truncated moments.
This savings also allows permutation testing and cross validation procedures for inferring the dimension of the eduction, which proved reasonably accurate in simulations. Finally, the application of the proposed methodology to socio-economic status (SES) index construction showed many advantages over common PCA-based indexes. 
In particular, the method not only helps get better predictions but also allows understanding the relations between the predictors and the response. 
More precisely, for the SES index, it gives varying weights capturing regional, historical and/or cultural differences, as well as various social measurement criteria (such as household per capita income or the poverty line), which it is not possible with PCA-based  methods. 
This property of the proposed method has relevant implications for the applied social analysis. 

Considering that many applications involve predictors of different natures (such as ordinal, continuous, and binary variables), further developments in SDR-based methods could be in this direction. 
The Principal Fitted Components (PFC) method for ordinal variables here proposed constitutes the first step to this extension. 
In particular,  the combination of ordinal and continuous predictors could be treated by taking all of them as continuous variables, where some of them are latent and the others are observable. 
Then, from the results here found, the reduction is identified, and the parameters can be estimated via maximum likelihood using the EM method on the latent variables and the PFC conventional method on the observed continuous variables. 
Nevertheless, the combination with binary variables requires a more exhaustive treatment, taking into account that the assumption of the existence of a latent normal variable on a binary variable may be naive and not make sense when the binary variable does not have a natural order (e.g., gender). 
Therefore, it is necessary to find a proper representation for the binary predictors, and search for a way to combine them with the other types of variables.

\newpage

\appendix

\section{The EM algorithm}\label{appendixA}

In order to simplify the notation we will always omit the conditioning on $\Omegak$ when taking expectations.
We will also omit the conditioning on some variables in the subscript. For instance, for any function $g$, we will call
\begin{align}\label{simple}
E_{\z_i} \big( g(\z_i) |y_i \big) &\doteq {E_{\z_i |y_i,{\tiny \Omegak}}} (g(\z_i) |y_i;\Omegak), \nonumber\\
E_{\z_i}  ( g(\z_i) |\x_i, y_i) &\doteq E_{\z_i| \x_i,y_i;\tiny\Omegak} (g(\z_i) |\x_i, y_i;\Omegak).
\end{align}

In order to obtain an explicit form of $Q$ we compute the conditional expectation of the joint log-likelihood. Following (\ref{simple}), we will write 
\[
Q(\Omegabf| \Omegak) = \sum_{i=1}^n E_{\z_i|y_i;\Omegak} \left[ \log f_{\x_i,\z_i}(\x_i,\z_i|y_i;\Omegabf)\big |y_i;\Omegak \right]  = \sum_{i=1}^n E_{\z_i} \left[ \log f_{\x_i,\z_i}(\x_i,\z_i|y_i;\Omegabf)\big |y_i\right].
\]

Therefore,
\small{
\begin{eqnarray} \label{Qnew}
&& Q(\Omegabf| \Omegak)  = \sum_{i=1}^n E_{\z_i} \left[ \log f_{\x_i,\z_i}(\x_i,\z_i|y_i;\Omegabf)\big |y_i \right] \nonumber \\ 
&&\;\;\;= \sum_{i=1}^n E_{\z_i}  \left[ \log \left((2\pi)^{-p/2} |\Deltabf|^{-1/2} e^{-\frac{1}{2}  \tr (\Deltabfs^{-1} 
(\z_i-\Deltabfs \alphabfs \xibfs  \fbar_{y_i})(\z_i-\Deltabfs\alphabfs \xibfs  
\fbar_{y_i})^T)}I_{\{\z_i \in C(\x_i,\Thetabfs)\}}\right)\Big|y_i \right] \nonumber  \\
&&\;\;\;= -\frac{pn}{2}\log(2\pi)- \frac{n}{2} \log |\Deltabf| \nonumber\\ &&\;\;\;- 
\frac{n}{2}\tr \left[ \Deltabf^{-1}\left(\frac{1}{n}\sum_{i=1}^n E_{\z_i}  (\z_i \z_i^T|\x_i, y_i) - \frac{2}{n} 
\Deltabf \alphabf \xibf  \sum_{i=1}^n \fbar_{y_i} E_{\z_i} (\z_i^T|\x_i,y_i)+ 
\frac{1}{n}\Deltabf\alphabf \xibf\sum_{i=1}^n  \fbar_{y_i}  f^T_{y_i} \xibf^T \alphabf^T \Deltabf 
\right) \right]  
\nonumber  \\&&\;\;\;=-\frac{pn}{2}\log(2\pi)- \frac{n}{2} \log |\Deltabf| - \frac{n}{2}\tr 
\left[\Deltabf^{-1} ( \Sbf - \frac{2 \Deltabf \alphabf \xibf  \F^T \M}{n} + \frac{\Deltabf\alphabf 
\xibf  \F^T \F \xibf^T \alphabf^T \Deltabf}{n})\right], 
\end{eqnarray}}
\noindent where $\Sbf \in \Rbb^{p\times p}$, $\F  \in \Rbb^{n\times r}$ and $\M  \in \Rbb^{n\times p}$ are given by $\Sbf = \frac{1}{n} \sum_{i=1}^n E_{\z_i}  (\z_i \z_i^T|\x_i, y_i)$, $\F^T = [\fbar_{y_1}, \ldots, \fbar_{y_n}]$ and $\M^T = [E_{\z_1} (\z_1|\x_1, y_1), \ldots, E_{\z_n} (\z_n|\x_n, y_n))]$, respectively.

\

\subsection{Maximizing the  Q-function (\ref{Qnew}).}

From (\ref{Qnew}), we have
\begin{align*}
Q(\parame) &=  -\frac{pn}{2}\log(2\pi)- \frac{n}{2} \log |\Deltabf| - \frac{n}{2}\tr 
(\Deltabf^{-1}\Sbf)+ \tr(\alphabf \xibf  \F^T \M)  - \frac{1}{2} \tr( \alphabf \xibf  \F^T \F \xibf^T 
\alphabf^T \Deltabf).
\end{align*}
Since $Q$ is a quadratic form in $\xibf$,  the maximum will be attained at
$\xibf^{(k)} = ( \alphabf^T  \Deltabf  \alphabf)^{-1} \alphabf^T \M^T \F 
(\F^T 
\F)^{-1}$. Replacing $\xibf^{(k)}$ in the Q-function, we obtain the partial log-likehood
\begin{align}\label{QQ}
Q(\Deltabf^{-1}, \alphabf) =& -\frac{pn}{2}\log(2\pi)- \frac{n}{2} \log |\Deltabf| - \frac{n}{2} 
\tr(\Deltabf^{-1} \Sbf) + \frac{1}{2} \tr 
\left[(\alphabf^T \Deltabf \alphabf)^{-1} \alphabf^T \M^T \F (\F^T \F)^{-1} \F^T \M \alphabf
\right] \nonumber \\ 
=& -\frac{pn}{2}\log(2\pi)- \frac{n}{2} \log |\Deltabf| - \frac{n}{2} 
\tr(\Deltabf^{-1} \Sbf) + \frac{n}{2} \tr 
\left[(\alphabf^T \Deltabf \alphabf)^{-1} \alphabf^T \Sbf_{\fit} \alphabf
\right],
\end{align}
where $\Sbf_{\fit} =  \frac{1}{n}\M^T \F (\F^T \F)^{-1} \F^T \M$. In order to maximize $Q$ with respect 
to $\Deltabf^{-1}$ observe that, by Proposition 5.14 in  \citep{Eaton1983}, if $\alphabf \in {\mathbb R}^{p\times d}$ is fixed and $\alphabf_0 \in {\mathbb R}^{p\times (p-d)} $ is the semi-orthogonal complement of
$\alphabf$, we have a one to one correspondence between 
$\Deltabf^{-1}$ and $(\Hbf_1,\Hbf_2,\Hbf_3)$, with
$\Hbf_1= \alphabf^{T}\Deltabf \alphabf$; $\Hbf_2=(\alphabf_0^{T} \Deltabf^{-1} \alphabf_0)^{-1}$ and
$\Hbf_3=(\alphabf^{T}\Deltabf \alphabf)^{-1} \alphabf^{T}\Deltabf \alphabf_0$.  From \citep{Eao73} and \citep{cook_forzani_2009} we have 
\begin{align}\label{esta}
\Deltabf^{-1}  &= \alphabf (\alphabf^T \Deltabf \alphabf)^{-1} \alphabf^T  + 
\Deltabf^{-1}  \alphabf_0 (\alphabf_0^T \Deltabf^{-1} \alphabf_0)^{-1}  \alphabf_0^T \Deltabf^{-1},\\
|\Deltabf | &= |\alphabf_0^T\Deltabf^{-1}  \alphabf_0|^{-1} |\alphabf^T \Deltabf \alphabf|.
\label{esta2}
\end{align}
Now the identity $(\alphabf^{T}\Deltabf \alphabf)^{-1} \alphabf^{T}\Deltabf \alphabf_0 = - \alphabf^{T}\Deltabf^{-1} \alphabf_0 (\alphabf_0^{T}\Deltabf^{-1} \alphabf_0)^{-1}$ implies that 
\[
\alphabf^T\Deltabf^{-1} \alphabf_0 = - (\alphabf^{T}\Deltabf \alphabf)^{-1} \alphabf^{T}\Deltabf \alphabf_0 \alphabf_0^{T} \Deltabf^{-1} \alphabf_0 =- \Hbf_3 \Hbf_2^{-1},
\]
which, together with $\alphabf \alphabf^T+\alphabf_0\alphabf_0^T = \mathbf I_d$, \[
\Deltabf^{-1} \alphabf_0 = \alphabf \alphabf^T\Deltabf^{-1} \alphabf_0+\alphabf_0\alphabf_0^T \Deltabf^{-1} \alphabf_0 = -\alphabf \Hbf_3 \Hbf_2^{-1} + \alphabf_0 \Hbf_2^{-1} = (-\alphabf \Hbf_3 + \alphabf_0)\Hbf_2^{-1}.
\]
With all this together in (\ref{esta}) we get
\begin{align}\label{esta1}
\Deltabf^{-1}  &= \alphabf \Hbf_1^{-1} \alphabf^T + (-\alphabf \Hbf_3 + \alphabf_0)\Hbf_2^{-1} (-\Hbf_3^T \alphabf^T+ \alphabf_0^T),
\end{align}
Therefore, finding $\widehat{\Hbf}_i$, $i=1,2,3$ is equivalent to finding $(\Deltabfi)^{(k)}$. 
In order to write the
Q-function in terms of $\Hbf_1$, $\Hbf_2$ and $\Hbf_3$, let us write $ \log 
|\Deltabf|$ and $\tr(\Deltabf^{-1} \Sbf)$ in terms of them. 
And therefore, using (\ref{esta1})  and (\ref{esta2}), the Q-function is then written in terms of $\Hbf_1$, $\Hbf_2$, and $\Hbf_3$ as 
\begin{align} \label{Qhaches}
Q(\Hbf_1, \Hbf_2,\Hbf_3, \alphabf) & = -\frac{pn}{2}\log(2\pi) - \frac{n}{2} \log |\Hbf_1| - 
\frac{n}{2} \log |\Hbf_2|    - \frac{n}{2} \tr \left[
\alphabf \Hbf_1^{-1} \alphabf^T (\Sbf -  \Sbf_{\fit}  ) \right] \nonumber \\ 
&\hspace{1cm} -  \frac{n}{2} \tr \left[( \alphabf_0 -  \alphabf \Hbf_3)\Hbf_2^{-1} ( 
\alphabf_0^T -  \Hbf_3^{T}  \alphabf^T)\Sbf \right].
\end{align}
Now, since $Q$ is quadratic in $\Hbf_3$, the maximum of $Q $ for $\Hbf_3$ is attained at
\begin{equation}\label{hache3}
\widehat{\Hbf}_3 = (\alphabf^{T} \Sbf  \alphabf)^{-1}(\alphabf^{T} \Sbf \alphabf_0).
\end{equation}
Replacing (\ref{hache3}) in (\ref{Qhaches}) and calling $\Sbf_{\res}= \Sbf - \Sbf_{\fit}$ (which is semidefinite positive), we have the partial log-likelihood function
\begin{align*} \label{Qhaches2}
Q(\Hbf_1, \Hbf_2, \alphabf) & = -\frac{pn}{2}\log(2\pi) - \frac{n}{2} \log |\Hbf_1|- 
\frac{n}{2} \log |\Hbf_2|    - \frac{n}{2} \tr \left[
\Hbf_1^{-1} \alphabf^T \Sbf_{\res} \alphabf \right] \nonumber \\ 
&\hspace{1cm} -  \frac{n}{2} \tr \left[\Hbf_2^{-1} ( 
\alphabf_0^T -  \alphabf_0^{T} \Sbf \alphabf(\alphabf^{T} \Sbf  \alphabf)^{-1} \alphabf^T) 
\Sbf ( \alphabf_0 -  \alphabf (\alphabf^{T} \Sbf  
\alphabf)^{-1}\alphabf^{T} \Sbf \alphabf_0)\right]. 
\end{align*}
The maximum of $Q$ over $\Hbf_1 $ and $\Hbf_2$ is attained at
\begin{eqnarray*}
\widehat{\Hbf}_1 &=& \alphabf^T \Sbf_{\res}\alphabf; \label{hache1}
\\\widehat{\Hbf}_2 &=&  \alphabf_0^{T} \Sbf \alphabf_0 
- \alphabf_0^{T} \Sbf \alphabf (\alphabf^{T} \Sbf \alphabf)^{-1} \alphabf^{T} \Sbf 
\alphabf_0 = (\alphabf_0^{T} \Sbf^{-1} \alphabf_0)^{-1}\label{hache2}.
\end{eqnarray*}
After substitution of the maximum for $\Hbf_1 $, $\Hbf_2$, and $\Hbf_3$ into (\ref{esta1}), we get that the maximum for $\Deltabfi$ is attained at
\begin{align*}
(\Deltabf^{-1})^{(k)}& =\Sbf^{-1}+ \alphabf (\alphabf^{T} \Sbf_{\res} \alphabf)^{-1} \alphabf^{T}  - \alphabf (\alphabf^{T} \Sbf \alphabf)^{-1} \alphabf^{T}. 
\end{align*}
Since this estimated matrix could not have unit elements in its diagonal, we scale it in order to have an unit-diagonal one. With this estimator of $\Deltabfi$, the partially maximized log-likelihood reads
\begin{eqnarray*} \label{Qhaches5}
Q(\alphabf) & = &
-\frac{pn}{2}\log(2\pi) - \frac{n}{2} \log 
|\alphabf^T \Sbf_{\res}\alphabf| -
\frac{n}{2} \log |(\alphabf_0^{T} \Sbf^{-1} \alphabf_0)^{-1}|- \frac{n 
d}{2} -   \frac{n (p-d)}{2}\nonumber
\\
&=& -\frac{pn}{2}[\log(2\pi) + 1] - \frac{n}{2} \log |\alphabf^T 
\Sbf_{\res}\alphabf|   - \frac{n}{2} \log |\Sbf|+ \frac{n}{2} \log |\alphabf^T \Sbf \alphabf|.
\end{eqnarray*}
where in the last equality we have used (\ref{esta2}). Finnally, the maximum in $\alphabf$ is attended at  
\begin{equation*}
\alphabf^{(k)} = \Sbf^{-1/2} \zetabfhat \Nbf,
\end{equation*}
where $\zetabfhat$ are the first $d$ eigenvectors of 
$ \Sbf^{-1/2}\Sbf_{\fit} \Sbf^{-1/2}$ and $\Nbf$ a matrix such that $\alphabfhat^T \alphabfhat = 
\mathbf I_d$. 

\newpage
\section{Approximating $\Sbf$ and $\M$}\label{appendiceC}
Given $(\Omegak,\Thetabf^{(k)})$ and $y_i$ fixed, we need to estimate $\Sbf$ and $\M$ in order to compute the Q-function.  
Each entry of matrix $\Sbf$ can be written as $s_{jk} = \sum_{i=1}^n 
E_{\z_i}(z_{i,j}z_{i,k}|\x_i)$ with $j,k=1,\ldots,p$. 
So, for $j=k$ we have the conditional second moment $E_{\z_i}(z_{i,j}^2|\x_i)$.
Following \citep{levina_2014}, when $j\neq k$ the terms $E_{\z_i} 
(z_{i,j}z_{i,k}|\x_i)$ can be approximated  by 
$E_{\z_i}(z_{i,j}z_{i,k}|\x_i)\approx E_{\z_i}(z_{i,j}|
\x_i)E_{\z_i}(z_{i,k}|\x_i)$. With this, we can obtain an estimator of $\Sbf$ through the 
estimation of first and second moments. 
The following is a modification of 
the procedure to compute these moments developed by \citep{levina_2014}, adapted to the case of conditional distributions.
We can write
$\x_i$ as $\x_i=(x_{i,j},\x_{i,-j})$ and $\z_i$ as $\z_i=(z_{i,j},\z_{i,-j})$ where 
$\x_{i,-j}=(x_{i,1},\ldots,x_{i,j-1},x_{i,j+1} ,\ldots, 
x_{i,p})$ and $\z_{i,-j}=(z_{i,1},\ldots,z_{i,j-1},z_{i,j+1} ,\ldots, z_{i,p})$. So, the first moment is   
\begin{eqnarray} \label{exp1}
E_{\z_i}(z_{i,j}|\x_i)&=&\int_{\Rbb^p} z_{i,j} f_{\z_i}(\z_i|\x_i) d\z_i \nonumber\\
&=& \int_{\Rbb^p} z_{i,j} f_{z_{i,j}}(z_{i,j}|\z_{i,-j},
x_{i,j}, \x_{i,j}) f_{\z_{i,-j}}(\z_{i,-j}| \x_i)d\z_i \nonumber\\
&=&\int_{\Rbb^{p-1}}\left[ \int_{\Rbb} z_{i,j}f_{z_{i,j}}(z_{i,j}|z_{i,-j},
x_{i,j})dz_{i,j}\right] f_{\z_{i,-j}}(\z_{i,-j}|\x_{i})d\z_{i,-j}\nonumber\\
&=&E_{\z_{i,-j}}\left\lbrace E_{z_{i,j}}(z_{i,j}|\z_{i,-j},x_{i,j})|\x_i \right\rbrace. 
\end{eqnarray}
In the same way,  the second moment can be written as
\begin{equation}\label{exp2}
\hspace{-1cm}E_{\z_i}(z^{2}_{i,j}|\x_i)=E_{\z_{i,-j}}\left\lbrace 
E_{z_{i,j}}(z^{2}_{i,j}|\z_{i,-j},x_{i,j})|\x_i \right\rbrace.
\end{equation}
Given $y_i$ and $\Omegak$, $\left(z_{i,1},z_{i,2} ,\ldots, z_{i,p}\right)$ has   
a multivariate normal distribution with mean $\mu_i= \Psikfy=\Deltak\alphabf^{\hspace{-0.05cm}(k-1)}\xibf^{(k-1)}$, and covariance matrix $\Deltak$.
Taking $\Deltak_{j,j}=1$, for each $j=1,\ldots,p$ we can write
\begin{center}
$\Deltak=\begin{pmatrix}
1 & \Deltak_{j,-j}\\
\Deltak_{-j,j} & \Deltak_{-j,-j}
\end{pmatrix}$\hspace*{1cm} and \hspace*{1cm}$\mu_i= \begin{pmatrix} (\Psibf^{(k-1)}\fbar_{y_i})_j\\
(\Psibf^{(k-1)}\fbar_{y_i})_{-j} \end{pmatrix}$,
\end{center}
and therefore the conditional distribution of $z_{i,j}$ given $\z_{i,-j}$ is
\begin{center}
$z_{i,j}|\z_{i,-j} \sim N(\widetilde{\mu}_{i,j},\widetilde{\delta}_{i,j})$ ,
\end{center}
where the mean is $\widetilde{\mu}_{i,j} = ( \Psikfy)_{j} + 
\Deltak_{j,-j}(\Deltak_{-j,-j})^{-1} \left(\z_{i,-j} - (\Psikfy)_{-j}\right)^{T}$ and the variance $\widetilde{\delta}^{2}_{i,j} = 1 - 
\Deltak_{j,-j}(\Deltak_{-j,-j})^{-1}\Deltak_{-j,j}$. In addition, the conditional distribution of $z_{i,j}$ on observed data $x_{i,j}$ is equivalent to conditioning on $z_{i,j}\in 
C({x_{i,j},\Thetabf})=[\theta^{(j)}_{x_{i,j}-1},\theta^{(j)}_{x_{i,j}})$, which follows a truncated 
normal distribution with density
\begin{align*}
f(z_{i,j}| C({x_{i,j},\Thetabf}), \z_{i,-j})= 
\frac{\frac{1}{\widetilde{\delta}_{i,j}}\phi(\widetilde{z}_{i,j})}{\Phi(\widetilde{\theta}^{(j)}_{x_{i,j}}) - 
\Phi(\widetilde{\theta}^{(j)}_{x_{i,j}-1})}I_{\left\lbrace z_{i,j} \in 
C({x_{i,j},\Thetabf})\right\rbrace}.
\end{align*}   
Here, $\widetilde{z}_{i,j}= (z_{i,j}-\widetilde{\mu}_{i,j}) / \widetilde{\delta}_{i,j}$,  
$\widetilde{\theta}^{(j)}_{x_{i,j}}=(\theta^{(j)}_{x_{i,j}} - \widetilde{\mu}_{i,j})/\widetilde{\delta}_{i,j}$ 
and $\widetilde{\theta}^{(j)}_{x_{i,j}-1}=(\theta^{(j)}_{x_{i,j}-1} - 
\widetilde{\mu}_{i,j})/\widetilde{\delta}_{i,j}$. From the moment generating function of the truncated normal distribution, the first and second moment of $z_{i,j}$ are given by
\begin{eqnarray} \label{exp3}
E\left(z_{i,j}|\z_{i,-j},x_{i,j}\right)=\widetilde{\mu}_{i,j}+\widetilde{\delta}_{i,j} a_{i,j} 
\label{exp31},\\
E\left(z^{2}_{i,j}|\z_{i,-j},x_{i,j}\right)=\widetilde{\mu}^{2}_{i,j}+\widetilde{\delta}^{2}_{i,j} + 2 
a_{i,j} \widetilde{\mu}_{i,j} \widetilde{\delta}_{i,j} 
+b_{i,j} \widetilde{\delta}^{2}_{i,j}\label{exp32},
\end{eqnarray}
where
\begin{center}
$a_{i,j}= 
\frac{\phi(\widetilde{\theta}^{(j)}_{x_{i,j}-1})-\phi(\widetilde{\theta}^{(j)}_{x_{i,j}})}{\Phi(\widetilde{
\theta}^{(j)}_{x_{i,j}})-\Phi(\widetilde{\theta}^{(j)}_{x_{i,j}-1})}$ ,\hspace*{1cm} 
$b_{i,j}=\frac{\widetilde{\theta}^{(j)}_{x_{i,j}-1}\phi(\widetilde{\theta}^{(j)}_{x_{i,j}-1})-\widetilde{\theta}
^{(j)}_{x_{i,j}}\phi(\widetilde{\theta}^{(j)}_{x_{i,j}})}{\Phi(\widetilde{\theta}^{(j)}_{x_{i,j}}
)-\Phi(\widetilde{\theta}^{(j)}_{x_{i,j}-1})}$.\\
\end{center}
Using (\ref{exp31}) and (\ref{exp32}) in (\ref{exp1}) and (\ref{exp2}), respectively, the first and 
second moments read
\begin{eqnarray}
E_{\z_i}(z_{i,j}|\x_i)=E_{\z_{i,-j}}\left(\widetilde{\mu}_{i,j}|\x_i \right)+ \widetilde{\delta}_{i,j} 
E_{\z_{i,-j}}\left( a_{i,j}|\x_{i}\right). \label{exp41}\\   
E_{\z_{i}}(z^{2}_{i,j}|\x_i)=E_{\z_{i,-j}}\left(\widetilde{\mu}^{2}_{i,j}|\x_i\right)+\widetilde{\delta}^{2}
_{i,j} + 2 \widetilde{\delta}_{i,j} E_{\z_{i,-j}}\left( a_{i,j}\widetilde{\mu}_{i,j}|\x_{i}\right)+ 
\widetilde{\delta}^{2}_{i,j} E_{\z_{i,-j}}\left( b_{i,j}|\x_{i}\right). \label{exp42}
\end{eqnarray}
Here $\widetilde{\mu}_{i,j}$ is linear in $\z_{i,-j}$ then, for $\Deltak$, $\alphabfhat$ and $y_i$ 
fixed, we have 
\begin{eqnarray}
\label{exp51}
& \hspace{-1cm} E_{\z_{i,-j}}\left(\widetilde{\mu}_{i,j}|\x_i \right)= (\Psikfy)_{j} + \Deltak_{j,-
j}(\Deltak_{-j,-j})^{-1} E_{\z_{i,-j}}\left(\z^{T}_{i,-j}|\x_{i}\right) - 
\Deltak_{j,-j}(\Deltak_{-j,-j})^{-1}
( \Psikfy)^{T}_{-j}\\
\nonumber
& \hspace{-1cm} E_{\z_{i,-j}}\left(\widetilde{\mu}^{2}_{i,j}|\x_i 
\right)=\left( \Psikfy\right)^{2}_{j} + 
2\left( \Psikfy\right)_{j}\Deltak_{j,-j}(\Deltak_{-j,-j})^{-1}\left[E_{\z_
{i,-j}}\left(\z_{i,-j}|
\x_i\right)-\left( \Psikfy\right)_{-j}\right]^{T}\\
\nonumber
& 
+\Deltak_{j,-j}(\Deltak_{-j,-j})^{-1}\Bigg[E_{\z_{i,-j}}\left(\z^{T}_{i,-j}\z_{i,-j}|\x_{i}
\right)- E_{\z_{i,-
j}}\left(\z^{T}_{i,-j}|\x_i\right)\left( \Psikfy\right)_{-j}\\
\label{exp52}
&-\left( \Psikfy\right)^{T}_{-j} E_{\z_{i,-j}}\left(\z_{i,-j}|
x_i\right)+\left( \Psikfy\right)^{T}_{-j}\left(\Psikfy\right)_{-j}\Bigg](\Deltak_{-j,-j})^{-1}(\Deltak_{-j,-j})^T
\end{eqnarray}

On the other hand, we have that the functions $a_{i,j}$ and $b_{i,j}$ are nonlinear in 
$\widetilde{\theta}^{j}_{x_{i,j}}$ and
$\widetilde{\theta}^{j}_{x_{i,j}-1}$ who are linear functions of $\widetilde{\mu}_{i,j}$ and thus of 
$\z_{i,-j}$. So, we can write $a_{i,j}$ and 
$b_{i,j}$ as $a_{i,j}(\z_{i,-j})$ and $b_{i,j}(\z_{i,-j})$. Conditioning on $\x_i$, $\z_{i,-j}$ has  
a truncated normal distribution with mean  
$\widetilde{\textbf{v}}_{i,-j} = E_{\z_{i,-j}}(\z_{i,-j}|\x_{i})$ and covariance matrix 
$\widetilde{\textbf{V}}=E_{\z_{i,-j}}\left((\z_{i,-
j}-\widetilde{\textbf{v}}_{i,-j})(\z_{i,-j}-\widetilde{\textbf{v}}_{i,-j})^{T} |\x_{i}\right)$. If  we 
assume that $a_{i,j}$ and $b_{i,j}$ have  continuous first partial derivatives, by the first order 
delta method we have that 
\begin{center}
$n^{1/2} \left\lbrace a_{i,j}(\z_{i,-j}|\x_{i}) - a_{i,j}(\widetilde{\textbf{v}}_{i,-j})\right\rbrace 
\xrightarrow{\mathcal{D}} N\left(0,\nabla a_{i,j}(\widetilde{\textbf{v}}_{i,-j}) 
\widetilde{\textbf{V}}\nabla^{T} a_{i,j}(\widetilde{\textbf{v}}_{i,-j})\right)$ and\\
$n^{1/2} \left\lbrace b_{i,j}(\z_{i,-j}|\x_{i}) -b_{i,j}( \widetilde{\textbf{v}}_{i,-j})\right\rbrace 
\xrightarrow{\mathcal{D}} N\left(0,\nabla b_{i,j}(\widetilde{\textbf{v}}_{i,-j}) 
\widetilde{\textbf{V}}\nabla^{T} b_{i,j}(\widetilde{\textbf{v}}_{i,-j})\right)$,
\end{center} 
so we can approximate the expectation $E_{\z_{i,-j}}\left( a_{i,j}|\x_{i}\right)$ with 
$a_{i,j}(\widetilde{\textbf{v}}_{i,-j})$ and $E_{\z_{i,-j}}\left( b_{i,j}|\x_{i}\right)$ with 
$b_{i,j}(\widetilde{\textbf{v}}_{i,-j})$, i.e. 
\begin{eqnarray} 
E_{\z_{i,-j}}\left( a_{i,j}|\x_{i}\right)\approx a_{i,j}(\widetilde{\textbf{v}}_{i,-
j})=\frac{\phi(\widetilde{\widetilde{\theta}}^{(j)}_{x_{i,j}-1})-\phi(\widetilde{\widetilde{\theta}}^{(j)}_{x_{i,j}}
)}
{\Phi(\widetilde{\widetilde{\theta}}^{(j)}_{x_{i,j}})-\Phi(\widetilde{\widetilde{\theta}}^{(j)}_{x_{i,j}-1})}\label{exp61}\\
E_{\z_{i,-j}}\left( b_{i,j}|\x_{i}\right)\approx b_{i,j}(\widetilde{\textbf{v}}_{i,-
j})=\frac{\widetilde{\widetilde{\theta}}^{(j)}_{x_{i,j}-1}\phi(\widetilde{\widetilde{\theta}}^{(j)}_{x_{i,j}-1}
)-\widetilde{\widetilde{\theta}}^{(j)}_{x_{i,j}}\phi(\widetilde{\widetilde{\theta}}^{(j)}_{x_{i,j}})}
{\Phi(\widetilde{\widetilde{\theta}}^{(j)}_{x_{i,j}})-\Phi(\widetilde{\widetilde{\theta}}^{(j)}_{x_{i,j}-1})}\label{exp62}
\end{eqnarray}
with
\begin{center}
$\widetilde{\widetilde{\theta}}^{(j)}_{x_{i,j}-k}=\frac{\widetilde{\theta}^{(j)}_{x_{i,j}-k} - 
E_{\z_{i,-j}}\left(\widetilde{\mu}_{i,j}|\x_i\right)} 
{\widetilde{\delta}_{i,j}}=\frac{\widetilde{\theta}^{(j)}_{x_{i,j}-k} - \left(\Psikfy\right)_j +\Deltak_{j,-j}(\Deltak_{-j,-j})^{-1}\left[\widetilde{\textbf{v}}_{i,-j}-\left(\Psikfy \right)_{-j}\right]^{T}}{\widetilde{\delta}_{i,j}}$
\end{center}
Using (\ref{exp51}), (\ref{exp52}), (\ref{exp61}), (\ref{exp62}) and the approximation 
$E_{\z_{i,-j}}\left( a_{i,j}\widetilde{\mu}_{i,j}|\x_{i}\right)\approx  E_{\z_{i,-j}}\left( 
a_{i,j}|\x_i\right) E_{\z_{i,-j}}\left(\widetilde{\mu}_{i,j}|\x_{i}\right)$, the conditional expectation in (\ref{exp41}) can be approximated by
\begin{align}\label{exp7}
E_{\z_{i}}\left(z_{i,j}|\x_i \right)&\approx ( \Psikfy)_{j} + 
\Deltak_{j,-j}(\Deltak_{-j,-j})^{-1} \Big[E_{\z_{i,-j}}\left(\z^{T}_{i,-j}|\x_{i}\right)-( \Psikfy)_{-j}\Big]^{T}\nonumber\\
&\hspace{1cm}+\widetilde{\delta}_{i,j} \frac{\phi(\widetilde{\widetilde{\theta}}^{(j)}_{x_{i,j}-1})-\phi(\widetilde{\widetilde{\theta}}^{(j)}_{x_{i,j}})}
{\Phi(\widetilde{\widetilde{\theta}}^{(j)}_{x_{i,j}})-\Phi(\widetilde{\widetilde{\theta}}^{(j)}_{x_{i,j}-1})},
\end{align}
and the second moment in (\ref{exp42}) by
\begin{align}\label{exp8}
&E_{\z_{i}}\left(z^{2}_{i,j}|\x_i \right) \approx \left( \Psikfy\right)^{2}_{j} + 
2\left( \Psikfy\right)_{j}\Deltak_{j,-j}(\Deltak_{-j,-j})^{-1}\left[E_{\z_
{i,-j}}\left(\z_{i,-j}|\x_i\right)-\left( \Psikfy\right)_{-j}\right]^{T} \nonumber\\
&\hspace{0.5cm}+\Deltak_{j,-j}(\Deltak_{-j,-j})^{-1}\Bigg[E_{\z_{i,-j}}\left(\z^{T}_{i,-j}\z_{i,-j}|\x_{i} \right)- E_{\z_{i,-j}}\left(\z^{T}_{i,-j}|\x_i\right)\left( \Psikfy\right)_{-j}\nonumber\\ &\hspace{0.5cm}-\left( \Psikfy\right)^{T}_{-j} E_{\z_{i,-j}}\left(\z_{i,-j}|x_i\right)+\left( \Psikfy\right)^{T}_{-j}\left( \Psikfy\right)_{-j}\Bigg](\Deltak_{-j,-j})^{-1}(\Deltak_{-j,-j})^T+\widetilde{\delta}_{i,j} \\
&\hspace{0.5cm}+ 2\widetilde{\delta}_{i,j}\frac{\phi(\widetilde{\widetilde{\theta}}^{(j)}_{x_{i,j}-1})-\phi(\widetilde{\widetilde{\theta}}^{(j)}_{x_{i,j}})}
{\Phi(\widetilde{\widetilde{\theta}}^{(j)}_{x_{i,j}})-\Phi(\widetilde{\widetilde{\theta}}^{(j)}_{x_{i,j}-1})} \Bigg[( \Psikfy)_{j} + \Deltak_{j,-j}(\Deltak_{-j,-j})^{-1} \Big[E_{\z_{i,-j}}\left(\z^{T}_{i,-j}|\x_{i}\right)-( \Psikfy)_{-j}\Big]^{T}\Bigg]\nonumber\\
&\hspace{0.5cm}+\widetilde{\delta}^{2}_{i,j}\frac{\widetilde{\widetilde{\theta}}^{(j)}_{x_{i,j}-1}\phi(\widetilde{\widetilde{\theta}}^{(j)}_{x_{i,j}-1})-\widetilde{\widetilde{\theta}}^{(j)}_{x_{i,j}}\phi(\widetilde{\widetilde{\theta}}^{(j)}_{x_{i,j}})}{\Phi(\widetilde{\widetilde{\theta}}^{(j)}_{x_{i,j}})-\Phi(\widetilde{\widetilde{\theta}}^{(j)}_{x_{i,j}-1})}. \nonumber
\end{align}
Equations (\ref{exp7}) and (\ref{exp8}) give recursive expressions for computing (iteratively) $\Sbf$ and $\M$, respectively. 

\newpage

\section{Description of Variables for SES index construction}\label{appendiceD}
The following variables are used to construct the SES indices:
\begin{itemize}
\item {\it Housing location}: indicates if the housing is located in a disadvantaged or vulnerable 
area. More precisely, it considers if housing: (i) is located in a shanty town, (ii) or/and near  to 
landfill sites,  (iii)  or/and  in a floodplain. It has 4 categories: 1 for houses 
that jointly present the  (i)-(iii) characteristics, 2 for housing presenting two of (i)-(iii), 3 
if  housing has only one of them, and 4 if the house has none of these characteristics.  
\item {\it Housing quality}: jointly contemplates the quality of roof, walls and 
floor  based on the CALMAT's methodology \citep{INDEC03} used in the population censuses of 
Argentina. It has 4 categories in increasing  order in terms of housing quality.
\item {\it Sources of cooking fuel}: indicates the kind of fuel used for cooking in the housing. 
It has 3 categories: 1 if the main source of cooking fuel in the housing is kerosene, wood or 
charcoal, 2 for bottled gas, and 3 for natural gas by pipeline. 
\item {\it Overcrowding}: characterizes the overcrowding by computing the ratio between rooms and 
number of household members. It has 4 categories: 1 if this ratio is less 
or equal than 1, 2 if the ratio is in the interval $(1,2]$, 3 if it is in $(2,3]$, and 4 if this ratio is 
greater than 3.
\item {\it Schooling}: indicates the formal education attained by the head of 
household. It has 7 categories: 1 if the head of household has no
formal education, 2 in the case of incomplete elementary level, 3 for complete 
elementary level, 4 for incomplete secondary school, 5 for a complete level of secondary school, 6 for an 
incomplete higher education and 7 if the head of household achieved a university or tertiary 
degree.
\item {\it Working hours}: describes the labor situation of head of household. It has 4 
categories: 1 for unemployment or inactive cases, 2 when the head of 
household works less than 40 hours per week, 3 for 40-45 per week working hours, and 4 when the 
head of household is employed for more than 45 hours per week. 
\item {\it Toilet drainage}: indicates the type of drainage of the housing. It has 4 categories: 
1 if drainage is a hole, 2 if drainage is only in a cesspool, 3 for cesspool and septic tank, and 4 for 
drain pipes in a public network.
\item {\it Toilet facility}: indicates the toilet facility available in the housing. It has 3 
categories: 1 for latrines, 2 for toilets without flush water, and 3 for flushing 
toilets.
\item {\it Toilet sharing}: indicates if the toilet is shared or not. It has 3 categories: 
1 if the toilet is shared with other housing, 2 if the toilet is  shared with 
other households into the same housing, and 3 if the toilet is used exclusively by the household.  
\item {\it Water location}: indicates the nearest location of drinking water. It 
has 3 categories: 1 if drinking water is gotten outside the plot of land of housing, 2 if 
water is into plot of land but outside of housing, and 3 of drinking water is obtained inside 
housing by pipe. 
\item {\it Water source}: indicates the source of the water in the housing. It has 3 categories: 
1 if drinking water comes from a hand pump or from a public tap shared with neighbours, 2 if 
drinking water is obtained by an automated drilling pump, and 3 for housing with piped drinking water.

\end{itemize}

\section*{Founding}
This work was supported by the SECTEI grant 2010-072-14, by the UNL grants 500-040, 501-499 and 500-062; by the CONICET grant PIP 742 and by the ANPCYT grant PICT 2012-2590.


\end{document}